\documentclass[12pt]{article}
\usepackage{fullpage}%varioref,multicol

\usepackage{bm}
\usepackage[dvipsnames]{xcolor}
\usepackage{amsmath,amssymb}
\usepackage{amsthm}
\usepackage{wrapfig}
\usepackage{epsfig}
\usepackage{graphicx}
\usepackage{makeidx}
\usepackage{multicol}
\usepackage{tikz}
\usepackage{algorithm}
\usepackage{algpseudocode}
\usepackage{soul}

\usepackage{array,multirow,booktabs,bbm} % for fancier tables

\usepackage{algorithmicx}
\usepackage{algpseudocode}
\usepackage{epstopdf}
\usepackage{mathrsfs,mathtools}
\usepackage{xcolor}
\usepackage{enumitem}
\newcommand{\N}{{\cal N}}

\newcommand{\D}{{\cal D}}

\newcommand{\R}{I\hspace{-1ex}R}

\graphicspath{{./Figs/}}
\newcommand{\bs}[1]{\boldsymbol{#1}}
\newcommand{\bmu}{\bs{\mu}}

\newcommand{\bc}{\bs{c}}

\newcommand{\bx}{\bs{x}}

\newcommand{\calP}{\mathcal{P}}
\newcommand{\calD}{\mathcal{D}}

\newcommand{\calN}{\mathcal{N}}

\DeclareMathOperator*{\argmax}{argmax}
\DeclareMathOperator*{\argmin}{argmin}

  \newtheorem{remark}{Remark}

\begin{document}
\graphicspath{{Figs/}}

\title{{L1-based reduced over collocation and hyper reduction for steady state and time-dependent nonlinear equations}}

\author{
Yanlai Chen\footnote{Department of Mathematics, University of Massachusetts Dartmouth, 285 Old Westport Road, North Dartmouth, MA 02747, USA. Email: {\tt{yanlai.chen@umassd.edu}}.}, \, 
Lijie Ji \footnote{School of Mathematical Sciences, Shanghai Jiao Tong University, Shanghai 200240, China. Email: {\tt sjtujidreamer@sjtu.edu.cn}. },\, 
Akil Narayan\footnote{Department of Mathematics, and Scientific Computing and Imaging (SCI) Institute, University of Utah, 72 S Central Campus Drive Room 3750, Salt Lake City, UT 84112 USA. Email: {\tt{akil@sci.utah.edu}}.},\,
Zhenli Xu \footnote{School of Mathematical Sciences, Institute of Natural Sciences,
and MOE-LSC,  Shanghai Jiao Tong University, Shanghai 200240, China. Email: {\tt xuzl@sjtu.edu.cn}. 
\newline ${   } \quad \ {  }$ L. Ji and Z. Xu acknowledge the support from grants NSFC 11571236 and 21773165 and HPC center of Shanghai Jiao Tong University. Y. Chen was partially supported by National Science Foundation grant DMS-1719698 and by AFOSR grant FA9550-18-1-0383. L. Ji is partly supported by China Scholarship Council (CSC, No.201906230067) during the author's one year visit at University of Massachusetts, Dartmouth. A.~Narayan was partially supported by NSF award DMS-1848508. This material is based upon work supported by the National Science Foundation under Grant No. DMS-1439786 and by the Simons Foundation Grant No. 50736 while Y.~Chen, L.~Ji, and A.~Narayan were in residence at the Institute for Computational and Experimental Research in Mathematics in Providence, RI, during the ``Model and dimension reduction in uncertain and dynamic systems" program.}
}

\date{\empty}

\maketitle

\begin{abstract}
The task of repeatedly solving parametrized partial differential equations (pPDEs) in, e.g. optimization, control, or interactive applications, makes it imperative to design highly efficient and equally accurate surrogate models. The reduced basis method (RBM) presents itself as such an option. Enabled by a mathematically rigorous error estimator, RBM carefully constructs a low-dimensional subspace of the parameter-induced high fidelity solution manifold from which an approximate solution is computed. It can improve efficiency by several orders of magnitudes leveraging an offline-online decomposition procedure. However, this decomposition, usually through the empirical interpolation method (EIM) when the PDE is nonlinear or its parameter dependence nonaffine, is either challenging to implement, or severely degrades online efficiency.

In this paper, we augment and extend the EIM approach as a direct solver, as opposed to an assistant, for solving nonlinear pPDEs on the reduced level. The resulting method, called Reduced Over-Collocation method (ROC), is stable and capable of avoiding the efficiency degradation inherent to a traditional application of EIM. Two critical ingredients of the scheme are collocation at about twice as many locations as the number of basis elements for the reduced approximation space, and an efficient L1-norm-based error indicator for the strategic selection of the parameter values to build the reduced solution space. Together, these two ingredients render the proposed L1-ROC scheme both offline- and online-efficient. A distinctive feature is that the efficiency degradation appearing in alternative RBM approaches that utilize EIM for nonlinear and nonaffine problems is circumvented, both in the offline and online stages. Numerical tests on different families of time-dependent and steady-state nonlinear problems demonstrate the high efficiency and accuracy of L1-ROC and its superior stability performance.
\end{abstract}

\section{Introduction}

Numerical simulations of systems, often parametrized, arising from various engineering and applied science disciplines are increasingly becoming of multi-query and/or real-time type.
For example, optimization and optimal control require multiple forward solves, interactive applications demand real-time responses. 
Design of fast numerical algorithms with certifiable accuracies for these settings has therefore continued to attract researchers' attention. The parameters delineating these systems may include boundary conditions, material properties, geometric settings, source properties etc. The wide variety, the complicated dependence of the system on these parameters, and their potential high dimensionality are the major challenges. In addition, the differential equations governing these equations may be nonlinear. 

The reduced basis method (RBM)  \cite{Quarteroni2015, HesthavenRozzaStammBook} has proved an effective option for this purpose. RBM was first introduced for nonlinear structure problem \cite{Almroth1978, noor1979reduced} in 1970s and has proven to be effective for linear evolution equations \cite{HaasdonkOhlberger}, viscous Burgers equation \cite{veroy2003reduced}, the Navier-Stokes equations \cite{deparis2009reduced}, and harmonic Maxwell's equation \cite{chen2010certified, chen2012certified}, just to name a few. The key to RBM's success in realizing orders-of-magnitude efficiency gain is an offline-online decomposition process where the basis selection and surrogate space construction are performed offline by a greedy algorithm, see review papers \cite{Rozza2008, Haasdonk2017Review} and monographs \cite{Quarteroni2015, HesthavenRozzaStammBook} for details. 
During the offline process, the necessary preparations for the online reduced solver are performed. The ultimate goal is that the complexity of the reduced solver, called upon in a potentially real-time fashion during the online stage, is independent of the number of degrees of freedom of the high-fidelity approximation of the basis functions.

\subsection{A key problem: Online efficiency degradation due to (D)EIM}

To achieve the efficiency goals of RBM, the {\em Empirical Interpolation Method} (EIM) or its discrete version (DEIM) \cite{Barrault2004, grepl2007efficient, ChaturantabutSorensen2010, PeherstorferButnaruWillcoxBungartz2014} is typically leveraged for nonaffine terms and/or nonlinear equations. However, EIM is often not feasible due to strong nonlinearity and/or nonaffinity of the problem. Even when it is feasible, performing (D)EIM can severely degrade the reduced solver's online efficiency when either the parameter dependence or the nonlinearity is complicated, such as when it encodes geometric variability \cite{chen2012certified, BenaceurEhrlacherErnMeunier2018}. The reason is that the online complexity is dependent on the number of terms resulting from the EIM decomposition.

Let us use a simple  system with a nonaffine parameter dependence as an example. Assume that we are solving a heat conduction problem with nonaffine parameter dependence $-\nabla \cdot \left( a(x; \bmu) \nabla u\right) = f$. For RBM to realize its intended efficiency gain, we would first apply EIM to approximate the function $a(x; \bmu)$ by a linear combination of $\bmu$-independent functions, $a(x; \bmu) \approx \sum_{q=1}^{Q_a}\theta_q(\bmu) a(x; \bmu^q)$, where $\{\bmu^q\}_{q=1}^{Q_a}$ is an ensemble, typically chosen through a greedy procedure. The number of terms $Q_a$ affects the online solver as follows, with the equation written in its weak form $a(u, v; \bmu) \coloneqq \left(a(x; \bmu) \nabla u, \nabla v\right) = (f, v)$. Once the offline learning stage identifies a reduced-order solution space spanned by e.g. full order solutions  $\{\xi^1, \dots, \xi^N\}$, the reduced solver is assembled for each $\bmu$, and the corresponding stiffness matrix has entries
\[
\left( a(\xi_i, \xi_j; \bmu) \right)_{i,j = 1}^N \coloneqq \left(a(x; \bmu) \nabla \xi_i, \nabla \xi_j\right)_{i,j = 1}^N = \sum_{q=1}^{Q_a} \theta_q(\bmu) \left(a(x; \bmu^q) \nabla \xi_i, \nabla \xi_j\right)_{i,j = 1}^N.
\]
The complexity of the online solver is therefore linearly dependent on the number of EIM terms $Q_a$, potentially suffering substantial reductions in efficiency compared to situations when EIM is not needed, i.e. $Q_a = 1$. The reason is that $Q_a$ can be prohibitively large (i.e. much larger than the reduced space dimension $N$) when the model involves geometric parametrization, see e.g. \cite{chen2012certified, BenaceurEhrlacherErnMeunier2018} even if the more efficient matrix version of EIM \cite{negri2015efficient} is adopted. As far as we are aware, efforts to mitigate this drawback are limited and underdeveloped.

\subsection{The proposed approach}

In this paper, we propose an L1-based reduced over-collocation method (L1-ROC) that is empirically stable and achieves full online efficiency without suffering $Q_a$-based efficiency degradation. Our main tools are an augmentation of EIM technique, a further leveraging of the collocation philosophy originally explored in \cite{ChenGottlieb2013}, and an extension of the L1 importance indicator proposed in \cite{JiangChenNarayan2019}. 
Let us summarize the two major ingredients of L1-ROC that, together, make the method able to circumvent this degradation. 

The first ingredient is a strategy to fully explore the EIM framework and adopt the collocation approach in contrast to variational approaches (i.e. Galerkin or Petrov-Galerkin) \cite{BennerGugercinWillcox2015, CarlbergBouMoslehFarhat2011, CarlbergBaroneAntil2017} when seeking the reduced solution. 
This so-called reduced collocation method is proposed and documented to work well in circumventing the EIM degradation for the reduced solver in previous work \cite{ChenGottlieb2013}. However its stability is lacking \cite{ChenGottliebMaday}. Our reduced over-collocation methods mitigate this stability defect by collocating at about twice as many locations as the dimension of the reduced order space. Half of these collocation points are identified from manipulation of a basis for this space: They interpolate the reduced solution (a linear combination of the basis elements). The other half are chosen according to a computational analysis of the reduced order residuals when these basis functions are identified during the offline procedure. They are present to ensure a good interpolation of the residual corresponding to an arbitrary parameter value when the reduced order space is used to solve the pPDE. This ingredient alone is not enough to achieve online and offline efficiency as the efficient calculation of the error estimators, critical for the construction of the reduced solution space, still relies on direct application of EIM. 

This challenge with computing error estimators is resolved by the second ingredient of the L1-ROC method, an efficient alternative for guiding the strategic selection of parameter values to build the reduced solution space. In particular, we utilize the recently introduced empirical L1 approach \cite{JiangChenNarayan2019} and extend it to time-dependent problems. Note that this approach is not a traditional rigorous error \textit{estimator}, and is instead simply an error \textit{indicator}. 

Together, these two ingredients render the L1-ROC scheme online-efficient (i.e. online cost is independent of the number of degrees of freedom of the high-fidelity truth approximation) and successfully circumvents the efficiency degradation of a direct EIM approach for nonlinear and nonaffine problems. Moreover, the L1-ROC method is highly efficient offline in that it requires minimal computation beyond the standard RBM cost of acquiring solution snapshots used to construct the reduced order space. As a consequence, the ``break-even'' number of simulations for the pPDE (minimum number of simulations that make the offline preparation stage worthwhile) is significantly smaller than traditional RBM. 
We test the algorithms on the viscous Burgers' equation \cite{veroy2003reduced} and various nonlinear convection diffusion reaction equations including the Poisson-Boltzmann equation. For all test problems, the L1-ROC is shown to have accuracy on par with the classical RBM while possessing much better efficiency due to the independence of the number of expansion terms resulting from the EIM decomposition. As examples, results for the steady-state and time-dependent cases of the diffusion with cubic reaction and the viscous Burgers' equation are shown.

\subsection{Other related techniques}

Popular model reduction techniques for linear time-dependent problems include Proper Orthogonal Decomposition (POD) \cite{Kunisch_Volkwein_POD}, system-theoretic approaches such as balanced truncation, moment matching or Hankel norm approximation \cite{BennerGugercinWillcox2015}. RBM stands out, for parametric problems in particular, with the availability of rigorous {\em a posteriori} error estimations, the resulting greedy algorithm, and the fact that it computes a number of full order solutions comparable to the theoretically smallest number, defined by the Kolmogorov $n$ width of the solution manifold.

The additional challenges posed by nonlinear problems are that a high dimensional reconstruction of the surrogate solution is usually needed for each evaluation of the nonlinearity. Sampling-based approximation techniques were developed to remedy this problem, including the Empirical Interpolation Method and its discrete variants \cite{Barrault2004, grepl2007efficient, ChaturantabutSorensen2010, PeherstorferButnaruWillcoxBungartz2014} and  Hyper-Reduction \cite{ryckelynck2005priori, Ryckelynck2009, CARLBERG2013623} which are known to be equivalent to DEIM under certain conditions \cite{Fritzen, Dimitriu2017}. Other approaches exist which include POD coupled with ``the best interpolation points'' approach  \cite{nguyen2008efficient, galbally2010non}, Gappy-POD \cite{Everson1995}, Missing Point Estimation (MPE) \cite{Astrid2008MissingPoint} or Gauss–Newton with approximated Tensors (GNAT) \cite{CarlbergBouMoslehFarhat2011, CARLBERG2013623}. Most of these methods work by first identifying a subset of the important features of the nonlinear function, and then constructing an approximation of the full solution based solely on an evaluation of these few components.

Our L1-ROC method can be viewed as adopting hyper reduction for {\em reduced} residual minimization. That is, instead of enforcing that the full residual is small in either a weak or strong formulation, one identifies its selected entries and ensures that an accurate evaluation of the residual on that subset is small. This is not the first time this type of idea is explored. For example, \cite{Astrid2008MissingPoint, astrid2004fast} uses a collocation of the original equations based on missing point interpolation and is followed by a Galerkin projection. The authors in \cite{ryckelynck2005priori} obtain the solution snapshots and collocation points through an adaptive algorithm in the finite element framework. It was also applied to nonlinear dynamical systems with randomly chosen collocation points \cite{bos2004accelerating}. 
However, the proposed L1-ROC differs from these existing works. The first distinctive feature is that the  basis functions and collocation points  are determined hierarchically via a greedy algorithm based on reduced residual minimization problems that gradually increase in size. The existing approaches obtain basis functions through POD-type techniques followed by computing collocation  points  all at once. The second distinctive feature is that the only step during the offline process that depends on the full order model is when we calculate a new high fidelity basis.

\vspace{3mm}

The paper is organized as follows. In Section \ref{sec:L1-ROC-Alg}, we introduce our L1-ROC method. Numerical results for two test problems, in both steady-state and time-dependent modes, are shown in Section \ref{num:final} to demonstrate the accuracy and efficiency of the scheme. Finally, concluding remarks are drawn in Section \ref{sec:conclusion}.

\section{The reduced over-collocation method}
\label{sec:L1-ROC-Alg}

In this section, we introduce the L1-ROC method for both steady state and time dependent problems.  
We first describe the problem we are solving. The framework of the online algorithm is then presented in Section \ref{sec:online}. Specification of part of the algorithm is postponed until the introduction of the ROC offline algorithm in Section \ref{sec:offline_c} which repeatedly calls the online solver to construct a surrogate solution space. 
The design of the main algorithm, the ROC approach, is detailed in Section \ref{sec:offline_oc}. To facilitate the reading of this and the following sections, we list our notation in Table \ref{tab:notation}.

We let $\mathcal{D} \subset \mathbb{R}^{p}$ be the domain for a $p$-dimensional parameter $\bmu$, and $\Omega \subset \mathbb{R}^{d}$ (for $d = 1 - \text{or} ~ 3$) be a bounded physical domain.  Given $\bm{\mu} \in \mathcal{D}$, and a Hilbert space $H$, the goal is to compute $u(\bm{\mu}) \coloneqq u(\bx;\bmu) \in H$ satisfying
\begin{equation}
\calP(u(\bx; \bmu);\bmu)-f(\bx)=0,~\bx \in \Omega,
\label{eq:steadymodel}
\end{equation}
or to compute time evolution of the transient problem
\begin{equation}
u_t +\calP(u(\bx; \bmu);\bmu)-f(\bx)=0,~\bx \in \Omega,
\label{eq:timemodel}
\end{equation}
with appropriate boundary (and initial) conditions. For example, for a stationary Laplace problem, the space $H$ is typically the Sobolev space $H^1(\Omega)$. Here, $\calP$ encodes a parametric second order partial differential operator that may include linear and nonlinear functions of $u(\bx; \bmu)$, $\nabla u(\bx; \bmu)$, and $\Delta u(\bx; \bmu)$. In the following, we will first focus on steady state problems \eqref{eq:steadymodel} and then extend the algorithm to time dependent case \eqref{eq:timemodel}.

To describe our algorithms, we first discretize the equation \eqref{eq:steadymodel} by a high-fidelity scheme (termed ``truth solver'' in the RB literature).  
In this paper, we adopt finite difference methods (FDM) for that purpose. However, extension to point-wise schemes such as spectral collocation is obvious, and to finite element methods is possible. We let $X^\N$ be a set of $\N$ collocation points on $\Omega$ at  which the equation is enforced on a discrete level. 
The discretized equation then becomes to find $u^\N(X^\N; \bmu)$, a discretization of the solution $u(\bmu)$ on the grid $X^\N$, such that 
\begin{equation}
\calP_{\N}(u^\N(X^\N; \bmu);\bmu)-f(X^\N)=0,
\label{eq:pdesystem}
\end{equation}
with $\nabla u(X^\N; \bmu)$, and $\Delta u(X^\N; \bmu)$ approximated by the numerical approximations $\nabla_h u(X^\N; \bmu)$, and $\Delta_h u(X^\N; \bmu)$. 
With a slight abuse of notation, we let $\N$ denote the number of the degrees of freedom in the solver,
even though the $\N$ points in $X^\N$ might include, e.g. points on a Dirichlet boundary that are not free.
\begin{table}[htbp]
  \begin{center}
  \resizebox{\textwidth}{!}{
    \renewcommand{\tabcolsep}{0.4cm}
    {\scriptsize
    \renewcommand{\arraystretch}{1.3}
    \renewcommand{\tabcolsep}{12pt}
    \begin{tabular}{@{}lp{0.8\textwidth}@{}}
      \toprule
      $\bmu = (\mu_1, \dots, \mu_p)$ & Parameter in $p$-dimensional parameter domain $\calD \subseteq \R^p$ \\
      $\Xi_{\rm{train}}$ & Parameter training set, a finite subset of $\mathcal{D}$ \\
      $u(\bmu)$ & Function-valued solution of a parameterized PDE on $\Omega \subset \mathbb{R}^{d}$\\
      $\calP(u(\bmu); \bmu)$ & A (nonlinear) PDE operator\\
 {$K$} & Number of finite difference intervals per direction of the physical domain\\
      $\mathcal{N} \approx K^d$ & Degrees of freedom (DoF) of a high-fidelity PDE discretization, the ``truth" solver \\ 
      $X^\calN$ & A size-$\calN$ (full) collocation grid\\
      $u^{\mathcal{N}}(\bmu)$ & Finite-dimensional truth solution\\
      $N$ & Number of reduced basis snapshots, $N \ll \mathcal{N}$\\
        $\bmu^j$ & ``Snapshot" parameter values, $j=1, \ldots, N$\\
      $\widehat{u}_n(\bmu)$ & Reduced basis solution in the $n$-dimensional RB space spanned by $\{u^{\mathcal{N}}(\bmu^1), \dots, u^{\mathcal{N}}(\bmu^n)\}$\\
      $e_n(\bmu)$ & Reduced basis solution error, equals $u^{\mathcal{N}}(\bmu) - \widehat{u}_n(\bmu)$ \\
      $\Delta_{{N}} \left(\bmu\right)$ & A residual-based error estimate (upper bound) for $\left\|e_N\left(\bmu\right)\right\|$ or an error/importance indicator\\
      $X^{N-1}_r=\{\bx^1_{**}, \dots, \bx^{N-1}_{**}\}$ & A size-$(N-1)$ reduced collocation grid, a subset of $X^\calN$ determined based on residuals\\
      $X^N_s = \{\bx^1_*, \dots, \bx^N_*\}$ & An additional size-$N$ reduced collocation grid, a subset of $X^\calN$ determined based on the solutions\\
      $X^M$ & A reduced collocation grid of size $M$ that is $X^{N-1}_r \cup X^N_s$\\
      $T$ & Final time for the time-dependent problems\\
      $\Delta t$ & Time stepsize for the time dependent problems\\
      $\calN_t$ & Total number of time levels, i.e. $\calN_t = \frac{T}{\Delta t}$\\
      $t^j$ & Time level $j$, $j=1, \ldots, \calN_t$\\
%      $t_i$ & Discretized time nodes, $i=1,2,\cdots,NT$.\\
      $\epsilon_{\mathrm{tol}}$ & Error estimate stopping tolerance in greedy sweep \\
      \midrule
      Offline component & The pre-computation phase, where the reduced solver is trained using a greedy selection of snapshots from the solution space\\
      Online component & The process of solving the offline-trained reduced problem, yielding the reduced order solution.\\
    \bottomrule
    \end{tabular}
  }
    }
  \end{center}
\caption{Notation and terminology used throughout this article.}\label{tab:notation}
\end{table}

\subsection{Online algorithm}

\label{sec:online}

The online component of the L1-ROC is essentially the same as the previously-introduced reduced collocation method \cite{ChenGottlieb2013} with the critical difference being that the number of collocation points is {\em larger} than the number of reduced basis snapshots. This {\em over-collocation} feature gives the method its name and provides additional stabilization of the online solver as we will observe in the numerical results. 

To describe the online algorithm, given $N$ selected parameters $\{\bmu^1, \dots, \bmu^N\}$, the corresponding {high fidelity truth approximations} $\{ u_n \coloneqq u^\N(X^\N;\bmu^n), 1 \le n \le N\}$, and $M \, (\ge N)$ collocation points formed from a subset of $X^\N$,
\begin{align*}
  X^M = \{\bx_*^1, \dots, \bx_*^M\}, \quad \mbox{ with }\bx_*^j  \mbox{ having index } i_j \mbox{ in }  X^\N, 
\end{align*}
we are able to perform the online algorithm, which we describe next. 
Note that, whenever there is no confusion, we are adopting the same notation for a function and its discrete representation in the form of a vector of its values at the grid points. These vectors $\{ u_n, 1 \le n \le N\}$ constitute the columns of basis matrix $W_{n} \in \mathbb{R}^{\N \times n}$ for $n \in \{1, \dots, N\}$. 
Furthermore, we denote the corresponding reduced representation of the basis space on the set $X^M$, by a matrix of the following form, 
\begin{align*}
W_{n,M} &= [u_{1}(X^M),  \ldots, u_{n}(X^M)] \in \mathbb{R}^{M \times n}, \quad \mbox{for } n = 1, \dots, N.\\
& =P_* W_{n},
\end{align*}
where the operator $P_* \in {\mathbb R}^{M \times \calN}$ is defined as,
\begin{align*}
P_* =\left[e_{i_1}, \cdots, e_{i_M}\right]^T,
\end{align*}
with $e_{i} \in \mathbb{R}^{\N \times 1}$ the $i$ th canonical unit vector.

Reduced approximations of the solution for any given parameter $\bmu$ are sought in the form,
$$\widehat{u}_n(\bmu) = W_{n} \bc_n (\bmu).$$ 
The condition for obtaining the coefficients $\bc_n (\bmu)$ is (a reduced version of) equation \eqref{eq:pdesystem} 
\begin{equation}
\calP_{\N}(W_{n} \bc_n (\bmu); \bmu) \approx f(X^\N).
\label{eq:reducedidea}
\end{equation}
Realizing that this is an over-determined system as we have in principle $n \ll \calN$, 
the authors of \cite{ChenGottlieb2013} proposed a Petrov Galerkin approach or collocation on $n$ points which produces a square system. 
The distinctive feature of what we propose in this paper  for locating the unknown coefficients $\bc_n(\bmu)$ is to minimizing  the residual of \eqref{eq:reducedidea} on the set of nodes $X^M$. Namely, we seek $\bc_n(\bmu)$ by solving the following optimization problem:
\begin{align}
\bc_n(\bmu) = \argmin_{\omega \in \mathbb{R}^n}\parallel P_* \left(\calP_{\N}(W_{n} \omega; \bmu)-f(X^\N)\right)\parallel_{\mathbb{R}^M}.
\label{pde:reduced}
\end{align}
We note that 
this is a nonlinear system of equations for $\bc_n$ with $\nabla_h \widehat{u}_n(\bmu)$ and $\Delta_h \widehat{u}_n(\bmu)$ computed on the full grid and then evaluated on the reduced grid $X^M$ according to
\begin{align*}
  \nabla_h\widehat{u}_n(\bmu) &=  P_* \left[\left(\nabla_h u_{1}\right),  \ldots, \left( \nabla_h u_{n}\right)\right] \bc_n(\bmu),\\
  \Delta_h\widehat{u}_n(\bmu)& =P_* \left[\left(\Delta_h u_{1}\right),  \ldots, \left( \Delta_h u_{n}\right)\right] \bc_n(\bmu).
\end{align*}
Iterative methods, such as Newton's method, will be used to solve for the coefficients $\bc_n(\bmu)$ online. 
The collocation nature of this scheme allows for solving this system with a cost only dependent on $M$ and $n$ even when $\calP_{\N}$ is nonlinear and nonaffine.  In particular, it is independent of the degrees of freedom $\N$ of the underlying truth solver. 
Indeed, the next section describes the offline procedure where the $N$ selected parameters $\{\bmu^1, \dots, \bmu^N\}$ are identified sequentially through a greedy algorithm. 
Once a $\bmu^j$ is determined,  we precompute as many quantities as possible so that minimal update is performed at each iteration of the iterative method. The online procedure of the nonlinear solve for obtaining $\bc_n(\bmu)$ from equation \eqref{pde:reduced} involves: 
\begin{itemize}
\item [1)] realizing/updating $W_{n,M}\bc_n$, $\nabla_h(W_{n,M})\bc_n$,  and $\Delta_h(W_{n,M})\bc_n$  at each iteration taking $O(M n)$ operations; 
\item [2)] calculating the forcing term $f(X^M)$ taking $O(M)$ operations; and 
\item [3)] solving the reduced linear systems at each iteration of the nonlinear solve taking $O(n^3)$ operations.
\end{itemize}

\subsection{Offline algorithm}

\label{sec:offline_c}

In this section, we describe the offline procedure of the reduced over collocation framework based on the L1-approach proposed in \cite{JiangChenNarayan2019}. 
The remaining ingredients of the offline procedure are identical with the traditional RBM algorithm \cite{Rozza2008, Haasdonk2017Review, Quarteroni2015, HesthavenRozzaStammBook}.

\subsubsection{A greedy algorithm based on an L1 importance indicator}
\label{sec:L1greedy}

We first briefly describe the procedure for  selecting the representative parameters $\bmu^1, \ldots, \bmu^N$ for constructing the solution space $W_N$. RBM utilizes a greedy scheme to iteratively construct $W_N$ relying on an  efficiently-computable error estimates that quantify the discrepancy between the dimension-$n$ RBM surrogate solution $\widehat{u}_n(\bmu)$ and the truth solution $u^\calN(\bmu)$.  Denoting such an estimate as $\Delta_n(\bmu)$, it traditionally satisfies $\Delta_n(\bmu) \geq \left\| \widehat{u}_n(\bmu) - u^\calN(\bmu)\right\|$. Assuming existence and computability of this error estimate, the greedy procedure for constructing $W_N$ then starts by selecting the first parameter $\bmu^1$ randomly from $\Xi_{\rm train}$ (a discretization of the parameter domain $\calD$) and obtaining
its corresponding high-fidelity truth approximation $u^\mathcal{N}(\bmu^1)$ to
form a (one-dimensional) RB space given by the range of $W_1 = \left[u^{\mathcal N}(\bmu^1)\right]$. Next, we
obtain an RB approximation $\widehat{u}_{n}(\boldsymbol{\mu})$ for each parameter in $\Xi_{\rm train}$ together with an error
bound $\Delta_n(\bmu)$. The greedy choice for the $(n+1)$th parameter $(n=1,\cdots,N-1)$ is made and the RB space augmented by
\begin{equation}
\label{eq:rbmgreedy}
  \bmu^{n+1} = \underset{\bm{\mu} \in \Xi_{\rm train}}{\argmax} \Delta_{n}(\bm{\mu}), \quad W_{n+1} = \left[ W_n \;\; u^\calN(\bmu^{n+1})\right] %W_n \oplus \{u^\calN(\bmu^{n+1})\}.
\end{equation}

The design and efficient implementation of the error bound $\Delta_n$ is usually accomplished with a residual-based {\em a posteriori} error estimate from the truth discretization. Mathematical rigor and implementational efficiency of this estimate are crucial for the accuracy of the reduced basis solution and its efficiency gain over the truth approximation. 
When $\calP(u; \bmu)$ is a linear operator, the Riesz representation theorem and a variational inequality imply that $\Delta_n$ can be taken as 
$$\Delta_n^R (\bmu) = \frac{\lVert f -
\mathcal{P}_\mathcal{N}(\widehat{u}_n; \bmu) \rVert_2} {\sqrt{\beta_{LB}(\bmu)}}, 
$$
which is a rigorous bound (with the $^R$-superscript denoting it is based on the full residual). Here $\beta_{LB}(\bmu)$ is a lower bound for the smallest eigenvalue of ${P}_\mathcal{N}(\bmu)^T {P}_\mathcal{N}(\bmu)$ with ${P}_\mathcal{N}(\bmu)$ being the matrix corresponding to the discretized linear operator $\mathcal{P}_\mathcal{N}(\cdot; \bmu)$.

Deriving the counterpart of this estimation for the general nonlinear equation is far from trivial. Moreover, even for linear equations, the robust evaluation of the residual norm in the numerator is delicate \cite{Casenave2014_M2AN, JiangChenNarayan2019}. We would also have to resort to an offline-online decomposition to retain efficiency which usually means application of EIM for nonlinear or nonaffine terms. This complication degrades, sometimes significantly \cite{BenaceurEhrlacherErnMeunier2018, negri2015efficient}, the online efficiency due to the large number of resulting EIM terms. 
What exacerbates the situation further is that the (parameter-dependent) stability factor $\beta_{LB}(\bmu)$ must be calculated by a computationally efficient procedure such as the successive constraint method \cite{HuynhSCM, HKCHP}. 
For these reasons, we are going to adopt the following empirical alternative, an {\em importance indicator} proposed in \cite{JiangChenNarayan2019}, in place of $\Delta_n^R$: 
$$\Delta_n^L(\bmu) = ||{\bf c}_n(\bmu)||_1.$$
The $^L$-superscript denotes that it is based on the L1-norm making our scheme L1-based. We note that this is not an error estimator because $\Delta_n^L$ does not decrease as we  increase $n$ since $\Delta_n^L(\bmu^i) = 1$ for $i \in \{1, \dots, n\}$. Nevertheless, we demonstrate that it is a reliable quantity to monitor when deciding which  representative parameters $\bmu^1, \dots, \bmu^N$ will form the surrogate space. We finish this subsection by pointing out that the calculation of $\Delta_n^L$ is independent of $\calN$ while naive approaches to evaluate the traditional estimator $\Delta_n^R$ for nonlinear problems would depend on $\calN$. This difference leads to the dramatic efficiency gain of the L1-ROC, as numerically confirmed in Section \ref{num:final}.
\begin{algorithm}[h]
\begin{algorithmic}[1]
\vspace{0.5ex}
\State Choose  $\bmu^1$ randomly from $\Xi_{\rm train}$, compute $u_1\coloneqq u^\N(X^\N;\bmu^1)$. 
\State Compute $\bx_*^1=\argmax_{\bx \in X^\N} |u_1|$, define $\xi_1 = u_1 / u_1(\bx_\ast^1)$. Let $i_1$ be the index of $\bx_*^1$ and $P_* = [e_{i_1}]^T$.
\State Initialize $m = n = 1, \, X^m = X^n_s =[\bx_*^1]$, $W_1 = \left\{\xi_1 \right\}, W_{1,m} = P_*W_1$, and $X_r^0 = \emptyset$. 
\State \mbox{\textbf{For}} $n = 2,\ldots, N$ 
\State $\quad\ \mbox{Solve }  \bc_{n-1} (\bmu) $ with $W_{n-1}, P_*$ and calculate $\Delta_{n-1}(\bmu)$  for all $\bmu \in \Xi_{\rm train}$.  
\State $\quad\ \mbox{Find } \bmu^{n} = \argmax_{\bmu \in \Xi_{\rm train}\backslash\left\{\bmu^i,i=1,\cdots,n-1\right\}}  \Delta_{n-1}(\bmu)$ and solve for $\xi_n \coloneqq u^\N(X^\N;\bmu^n)$. 
\State $\quad\ \mbox{Compute an interpolatory residual for } \xi_n: \, \mbox{find } \{\alpha_j\} \mbox{ and let } \xi_n = \xi_n - \sum_{j=1}^{n-1}\alpha_j \xi_j$ \mbox{so that } $\xi_n(X^{n-1}_s)=0$.
\State $\quad\ \mbox{Find }$ $\bx_*^n=\argmax_{\bx \in X^\N/X^m} |\xi_n(\bx)|$, $\xi_n=\xi_n/ \xi_n(\bx_*^n)$, and let $X^n_s = X^{n-1}_s \cup \{\bx_*^n\}$, and $i_1$ be the index of $\bx_*^n$.
\State $\quad\ \mbox{Form the full residual vector } r_{n-1} =\calP_{\N}(\widehat{u}_{n-1}(\bmu^n);\bmu^n)-f(X^\N)$ and compute its interpolatory residual: $\, \mbox{find } \{\alpha_j\} \mbox{ and let } r_{n-1} = r_{n-1} - \sum_{j=1}^{n-2}\alpha_j r_j$ \mbox{so that } $r_{n-1}(X_r^{n-2})=0$. Find $\bx^{n-1}_{**}=\argmax_{\bx \in X^\N/\left\{X^m,\bx_*^n\right\}} |r_{n-1}(\bx)|$. Let $r_{n-1}=r_{n-1}/ r_{n-1}(\bx^{n-1}_{**})$, and $X^{n-1}_r =X^{n-2}_r \cup \{\bx^{n-1}_{**}\}$ and $i_2$ is the index of $\bx_{**}^{n-1}$.
\State $\quad\ \mbox{Update } W_{ n} = \{W_{n-1}, \xi_n\}, m=2n- 1, X^m = X^n_s \cup X_r^{n-1}, P_*= P_*\cup [e_{i_1},e_{i_2}]^T$.
\State \mbox{\textbf{End For}}
\end{algorithmic}
\caption{Offline: construction of $W_N$ and the  collocation set $X^{2N - 1} = X^N_s \cup X^{N-1}_r$.} 
\label{alg:c:plus:offline2}
\end{algorithm}

\subsubsection{Construction of the reduced over-collocation set $X^M$}

\label{sec:offline_oc}

Let us now describe how we determine the reduced collocation set $X^M$ to complete the offline algorithm. 
Toward that end, we first describe the construction of two sets. The first one, denoted by $X^N_s$, 
consists of the maximizers from the EIM procedure used on the orthonormalized columns of $W_N$, which are computed as pivots from an LU decomposition. Realizing the importance of controlling the residuals when solving equations, we need to represent the residuals well on the reduced grid. For that purpose, we introduce a second set of points, examines the residual of the RB solution at the chosen $\bmu^{n}$ when only $n-1$ basis elements are used, 
\begin{equation}
r_{n-1} =\calP_{\N}(\widehat{u}_{n-1}(\bmu^n);\bmu^n)-f(X^\N), \quad n \in \{2, \dots, N\}.
\label{eq:offlineresidual}
\end{equation}
We next take these $N-1$ residual vectors and perform 
an EIM procedure on them. 
The $N-1$ maximizers from this procedure form the second set which is denoted $X^{N-1}_r$. 
The reduced collocation approach in \cite{ChenGottlieb2013} is a specialization that takes $M = N$ and $X^M = X^N_s$. The resulting $M=N$ reduced scheme can be unstable particularly when high accuracy (i.e. large $N$) of the reduced solution is desired. It can be resolved in special cases by an analytical preconditioning approach \cite{ChenGottliebMaday}. 
The second obvious choice of $X^M$ is to append $X^{N-1}_r$ with one more point such as the maximizer of the first basis. Numerical tests (not reported in this paper) also reveal instability of this scheme. 

The stabilization mechanism and name of the reduced over-collocation methods, outlined in Algorithm \ref{alg:c:plus:offline2}, come from the fact that we combine these two choices by taking $$M = 2N - 1 \mbox{ and } X^M = X^N_s \cup X^{N-1}_r,$$ and solving a least squares problem on the reduced level by collocating on about twice as many points as the number of basis in the RB space. Note that the first basis function has no accompanying residual vector \eqref{eq:offlineresidual}, so that from the second onward there are two collocation points selected whenever a new parameter is identified by the greedy algorithm.

\subsection{Extension of L1-ROC for time dependent problems}
\label{sec:offline:time}

For the time-dependent problem \eqref{eq:timemodel}, the semi-discretized L1-ROC solver remains identical to the steady-state case. That is, we seek the reduced approximation of the solution for any given parameter $\bmu$  in the form of
$$\widehat{u}_n(\bmu, t) = W_{n} \bc_n (\bmu, t).$$ 
The unknown coefficients $\bc_n(\bmu, t) \in {\mathbb R}^{n \times 1}$ is obtained by solving the following optimization problem:
\begin{align}
\bc_n(\bmu, t) = \argmin_{\omega \in \mathbb{R}^n}\parallel P_* \left(W_{n} \omega + \calP_{\N}(W_{n} \omega; \bmu)-f(X^\N)\right)\parallel_{\mathbb{R}^M}.
\label{pde:reduced_t}
\end{align}

To discretize in time, our L1-ROC aligns with the parameter-time greedy framework \cite{grepl2007efficient, grepl2005posteriori}, as opposed to POD \cite{Kunisch_Volkwein_POD, nguyen2008efficient} or POD-greedy \cite{grepl2012}. We discretize the time and denotes the (full) set of temporal nodes as ${\mathcal T}_f \coloneqq \{t_i: \, i =0, \cdots, \calN_t\}$ with $t_0$ being the initial time and $\calN_t =T/ \Delta t$ where $\Delta t$ is the temporal step-size. We extend the L1-based importance indicator of \cite{JiangChenNarayan2019} to the time-dependent case here. Toward that end, we define a {\em reduced} set of temporal nodes ${\mathcal T}_r$ that starts from the empty set and is gradually enriched in the greedy algorithm. 

To initiate the reduced solver construction we start with a deterministically or randomly chosen $\bmu^1$ (similar to the steady-state case) and invoke the truth solver to obtain the snapshots $\{u^\calN(t_i, x; \bmu^1)\}_{i=0}^{\calN_t}$. ${\mathcal T}_r$ is initiated by the time instant when the corresponding snapshot has the largest variation. That is,
\[
  {\mathcal T}_r = \{t_{\bmu^1}^{1}\} \mbox{ where } t_{\bmu^1}^{1} = \argmax_{t \in {\mathcal T}_f} \left( \max_{x \in X^\N} u^\calN(t, x; \bmu^1)- \min_{x \in X^\N} u^\calN(t, x; \bmu^1) \right).
\]
The RB space $W_1$ is initiated with $u^\calN(t_{\bmu^1}^{1}, x; \bmu^1)$. The (first) collocation point is set to be the EIM point of this first basis, i.e. the spatial maximizer of $|u^\calN(t_{\bmu^1}^{1}, x; \bmu^1)|$,
\[
  \bx_\ast^1
= \argmax_{x \in X^\N} |u^\calN(t_{\bmu^1}^{1}, x; \bmu^1)|.
\]

Once these ingredients are in place after the first pair $(\bmu^1, t^1_{\bmu^1})$ is determined, we can solve the reduced problem \eqref{pde:reduced_t} for every $\bmu \in \Xi _{\rm train}$.  Similar to the traditional greedy algorithm, the next step is to determine the subsequent $(\bmu, t)$ pairs. Our greedy algorithm manifests itself in the following three aspects:
\begin{itemize}
\item {\bf Greedy in $\bmu$:} 
We define the following importance indicator for each $\bmu$ after its corresponding (reduced) solver of \eqref{pde:reduced_t} is performed,
\begin{align}
\Delta_n^{Lt}(\bmu) \coloneqq \max_{t \in {\mathcal T}_r} \{{\lVert\bc_{n}(\bmu, t)\rVert}_1 \}.
\label{eq:delta_t}
\end{align}

We note that: 1) the maximization is done only on the reduced temporal grid ${\mathcal T}_r$ which is much smaller than the full temporal grid ${\mathcal T}_f$; 2) the signature feature of the L1-based approach carries over to the time-dependent case in that the indicator requires nothing more than the reduced solution coefficients. Our greedy choice for the $\bmu$-component of the $(\bmu, t)$ pair is through maximizing $\Delta_n^{Lt}(\bmu)$ over the training set $\Xi _{\rm train}$:
\[
  \bmu^{n+1} = \argmax_{\bmu \in \Xi _{\rm train}} \Delta_n^{Lt}(\bmu).
\]
\item {\bf Greedy in $t$:} Next, the $t$-component of the $(\bmu, t)$ pair is determined and the set ${\mathcal T}_r$ enriched with a new temporal node through a greedy choice as well.  Given the greedy choice $\bmu^{n+1}$ and the {\em reduced} solution $\widehat{u}_n(\bmu^{n+1}, t) = W_{n} \bc_n (\bmu^{n+1}, t)$ for all time levels $t \in {\mathcal T}_f$, we compute the {\em full} residual vectors $r_n \in \mathbb{R}^{\N \times 1}$ for this $\bmu^{n+1}$. 
The greedy $t$-choice is given by 
\begin{align}
t_{\bmu^{n+1}}^{k_{\bmu^{n+1}}} \coloneqq \argmax_{t \in {\mathcal T}_f}\left\{\varepsilon( t;\bmu) \coloneqq {\lVert r_n (t;\bmu^{n+1})\lVert}_\infty\right\}, \mbox{ and } {\mathcal T}_r = {\mathcal T}_r \bigcup \{t_{\bmu^{n+1}}^{k_{\bmu^{n+1}}}\}.
\label{eq:reducedresidual}
\end{align}
Here, $k_{\bmu^{n+1}} \ge 1$ is introduced to account for the possibility that multiple temporal nodes might be selected for the same $\bmu$, at different rounds of the greedy algorithm. 
We note in particular that, consistent with typical greedy scheme, we choose one (as opposed to multiple) maximizer in \eqref{eq:reducedresidual}. However, as we proceed with building up the reduced solution space, the same $\bmu$ (and a different temporal node) may be chosen by the greedy algorithm at a later step due to the lack of resolution of its corresponding temporal history. 
\item {\bf $X^M$ expansion:} Once a new greedy pair $(\bmu^{n+1}, t_{\bmu^{n+1}}^{k_{\bmu^{n+1}}})$ is fixed, we solve for the truth approximations $u(t, X^\N; \bmu^{n+1})$ for $t \le t_{\bmu^{n+1}}^{k_{\bmu^{n+1}}}$. The expansion of $X^M$ by two more colocation points, with one from the EIM procedure of the solution $u(t_{\bmu^{n+1}}^{k_{\bmu^{n+1}}}, X^\N; \bmu^{n+1})$ and the other from that of the residual $r_n(t_{\bmu^{n+1}}^{k_{\bmu^{n+1}}}; \bmu^{n+1})$, is identical to the steady state case. 
\end{itemize}

The full offline algorithm is seen in Algorithm.\ref{alg:c:plus:offline:time}.
\begin{algorithm}[h]
\begin{algorithmic}[1]
\vspace{0.5ex}
\State Choose $\bmu^1$, and set $k_{\bmu^1} = 1$ the first temporal node to be $t_{\bmu^1}^{k_{\bmu^1}} = \argmax_{t \in {\mathcal T}_f} \left( \max_x u^\calN(t, x; \bmu^1)- \min_x u^\calN(t, x; \bmu^1) \right)$. Define $\xi_1\coloneqq u^\mathcal{N}(t_{\bmu^1}^{k_{\bmu^1}}, X^\N;\bmu^1)$. 
\State Find $\bx_*^1=\argmax_{\bx \in X^\N} |\xi_1|$, and let $P_* = [e_{i_1}]^T$, where $i_1$ is the index of $\bx_*^1$.
\State Initialize $m = n = 1, \, X^m = X^n_s =\{\bx_*^1\}$, $W_1 = \left\{\xi_1 \right\}, W_{1,m} = P_*W_1$, and $X_r^0 = \emptyset$.  %\\[0.5ex]
\State  \mbox{\textbf{For}} $n = 2,\ldots, N$  %\\%[0.5ex]
\State  \quad Solve   the reduced problem for $\bc_{n-1} (\bmu, t_k) $. 
\State  $\quad\ \mbox{Find } \bmu^n = \argmax_{\bmu \in \Xi_{\rm train}}  \Delta_{n-1}^{Lt}(\bmu)$, and a new temporal node $t_{\bmu^n}^{k_{\bmu^n}} = \arg \max_{t \in {\mathcal T}_f} {\varepsilon(t;\bmu^n)}$. %\\%[.5ex]
\State  $\quad\ \mbox{Solve }  \xi_n=u^\mathcal{N} (t_{\bmu^n}^{k_{\bmu^n}}, X^\N;\bmu^n)$. %\\%[0.5ex]
\State  $\quad\ \mbox{Compute an interpolatory residual for } \xi_n: \, \mbox{find } \{\alpha_j\} \mbox{ and let } \xi_n = \xi_n - \sum_{j=1}^{n-1}\alpha_j \xi_j$ \mbox{so that } $\xi_n(X^{n-1}_s)=0$. Find $\bx_*^n=\argmax_{\bx \in X^\N/X^m} |\xi_n|$, $\xi_n=\xi_n/ \xi_n(\bx_*^n)$. Let $X^n_s = X^{n-1}_s \cup \{\bx_*^n\}$, and $i_1$ be the index of $\bx_*^n$.%\\[0.5ex]
\State  $\quad\ \mbox{Form the full residual vector } r_{n-1} = \left(\widehat{u}_{n-1}\right)_t(t_{\bmu^n}^{k_{\bmu^n}};\bmu^n) + \calP_{\N}(X^\N, \widehat{u}_{n-1}(t_{\bmu^n}^{k_{\bmu^n}};\bmu^n);\bmu^n)-f(X^\N, t_{\bmu^n}^{k_{\bmu^n}})$. 
Compute an interpolatory residual $r_{n-1}: \mbox{find } \{\alpha_j\} \mbox{ and let } r_{n-1} = r_{n-1} - \sum_{j=1}^{n-2}\alpha_j r_j$ \mbox{so that } $r_{n-1}(X^{n-2}_r)=0$. Find $\bx_{**}^n=\argmax_{\bx \in X^\N/\left\{X^m,\bx_{*}^n\right\}} |r_{n-1}|$.Let  $r_{n-1}=r_{n-1}/ r_{n-1}(\bx_{**}^n)$, and $X^{n-1}_r =X^{n-2}_r \cup \{\bx_{**}^n\}$. $i_2$ is the index of $\bx_{**}^n$.%\\%[0.5ex]
\State  $\quad\ \mbox{Update } W_{n} = \{W_{n-1},\xi_n\}, m=2n- 1, X^m = X^n_s \cup X_r^{n-1},P_* = [P_*; \,\, (e_{i_1})^T; \, \, (e_{i_2})^T]$. %\\%[0.5ex]
\State  \mbox{\textbf{End For}}
\end{algorithmic}
\caption{L1-ROC algorithm for time dependent problems}
\label{alg:c:plus:offline:time}
\end{algorithm}

\section{Numerical results}
\label{num:final}

In this section, we present the numerical results of the L1-ROC method applied to the nonlinear steady-state and time-dependent problems, in Sections \ref{numerics:steady} and \ref{numerics:timedep} respectively. The equations we test, in each section, include the classical viscous Burgers' equation and nonlinear convection diffusion reaction equations. 

\subsection{L1-ROC for steady-state nonlinear problems}
\label{numerics:steady}

\subsubsection{Viscous Burgers' equation}
First, we show the results of our algorithm applied to the one-dimensional (viscous) Burgers' equation,
  \begin{equation}
 \begin{split}
 u u_x & = \bmu u_{xx},\\
 u(x=-1) & = 1, \,\,\,\, u(x=1) = -1.
 \end{split}
 \label{eq:burgers}
\end{equation} 
Here the viscosity parameter $\bmu$ varies on the interval $\calD = [0.05,1]$. The computational domain $[-1,1]$ is divided uniformly into $\calN+1$ intervals with grid points denoted by 
\[
\{x_0, x_1, \dots, x_{\calN+1}\}.
\]
With $h = \frac{2}{\calN+1}$, the following  finite difference discretization based on the conservative form of equation \eqref{eq:burgers}, $\left( \frac{u^2}{2} \right)_x - \bmu u_{xx} = 0$, is then used 
\begin{equation}
\frac{u_{i+1}^2 -u_{i-1}^2}{4 h}- \bmu \frac{u_{i-1} -2u_i +u_{i+1}}{h^2}=0, \quad i \in \{1, \dots, \calN\}. 
\label{eq:conservative}
\end{equation}
This leads to a nonlinear truth solver of size $\calN$ to resolve \eqref{eq:burgers}.
The parameter domain $\calD$ is sampled 50 times logarithmically spaced to form the training set for the Offline procedure. 

We test our method on a subset of $\Xi_{\rm test}$ of $\calD$ that has empty intersection with the training set $\Xi_{\rm train}$. 
We compute the relative errors $E(n)$ over all $\bmu$ in $\Xi_{\rm test}$ of the reduced basis solution using $n$ basis functions, $\widehat{u}_n(\bmu)$, in comparison to the high fidelity truth approximation. That is,
\begin{equation}
E(n) =  \max_{\bmu \in \Xi_{\rm test}}\left\{\frac{\| u(\bmu) - \widehat{u}_n(\bmu)\|_\infty}{\|u\|_{L^\infty(\Xi_{\rm test}, L^\infty(\Omega))}}\right\}
\label{eq:error:steadyburger}
\end{equation}
where 
\[
||u||_{L^\infty(\Xi_{\rm test}, L^\infty(\Omega))} = \max_{\bmu \in \Xi_{\rm test}} \|u(\bmu) \|_\infty.
\]

\begin{figure}[!htb]
\centering
\includegraphics[height=0.19\textheight]{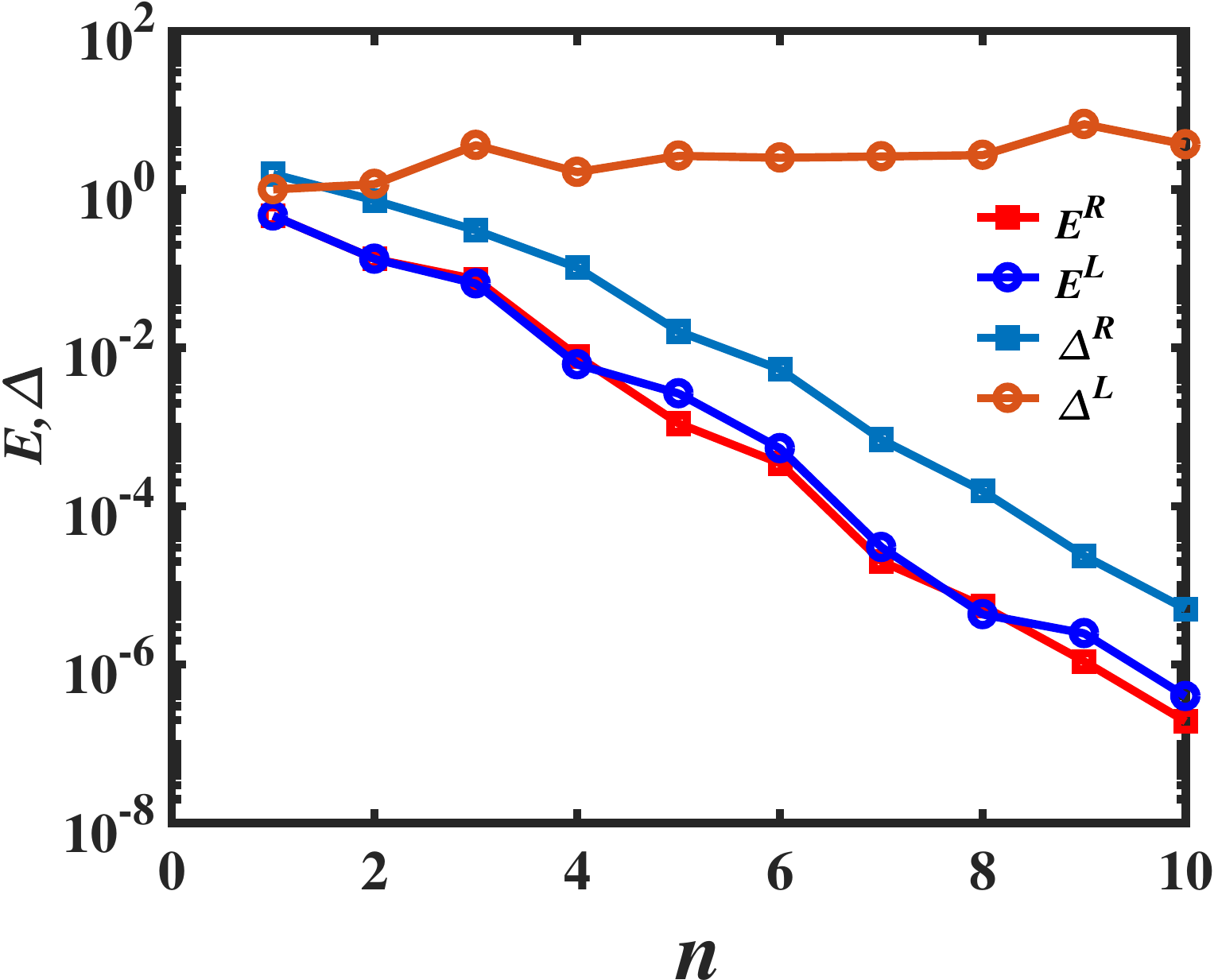}
\includegraphics[height=0.19\textheight]{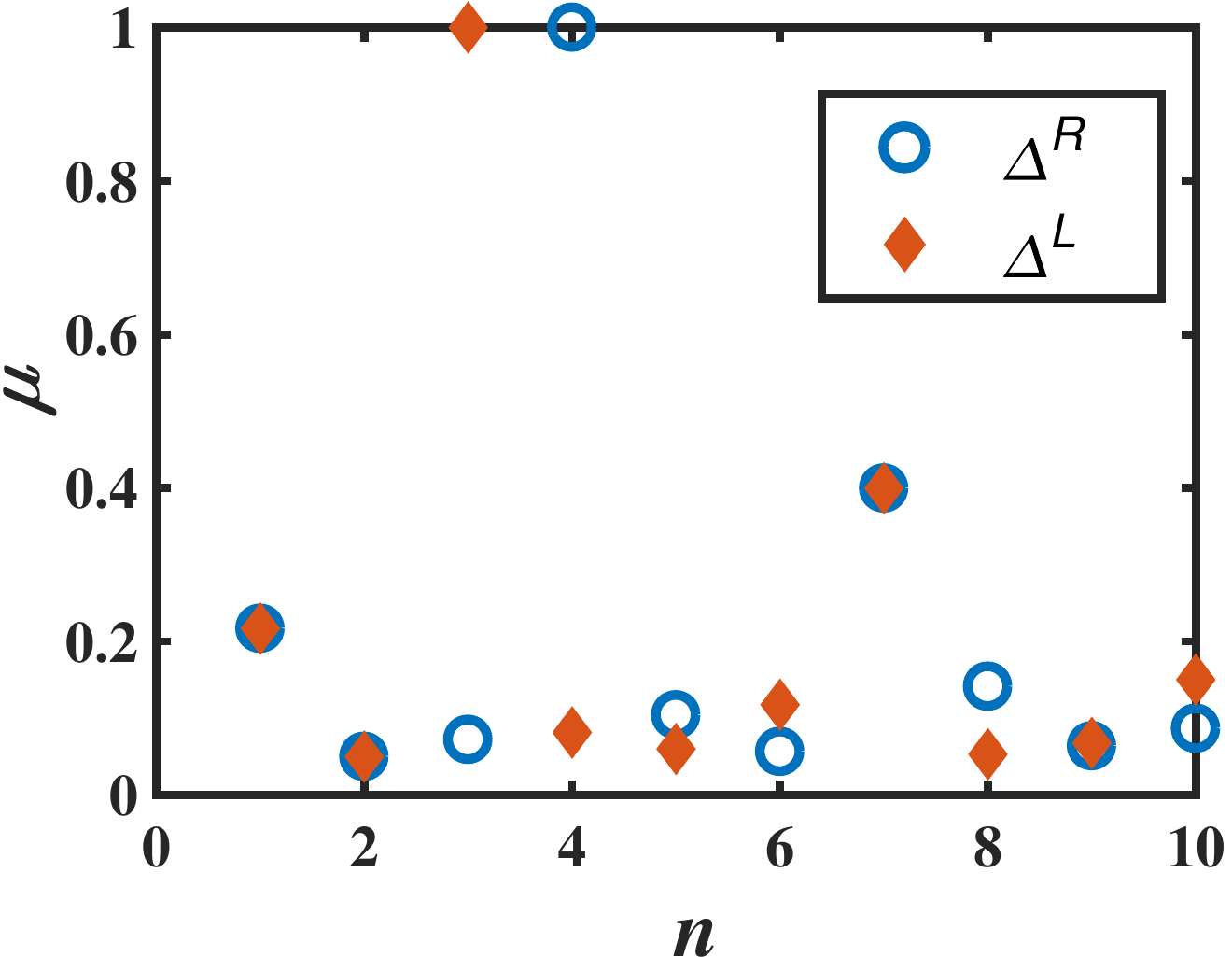}
\includegraphics[height=0.19\textheight]{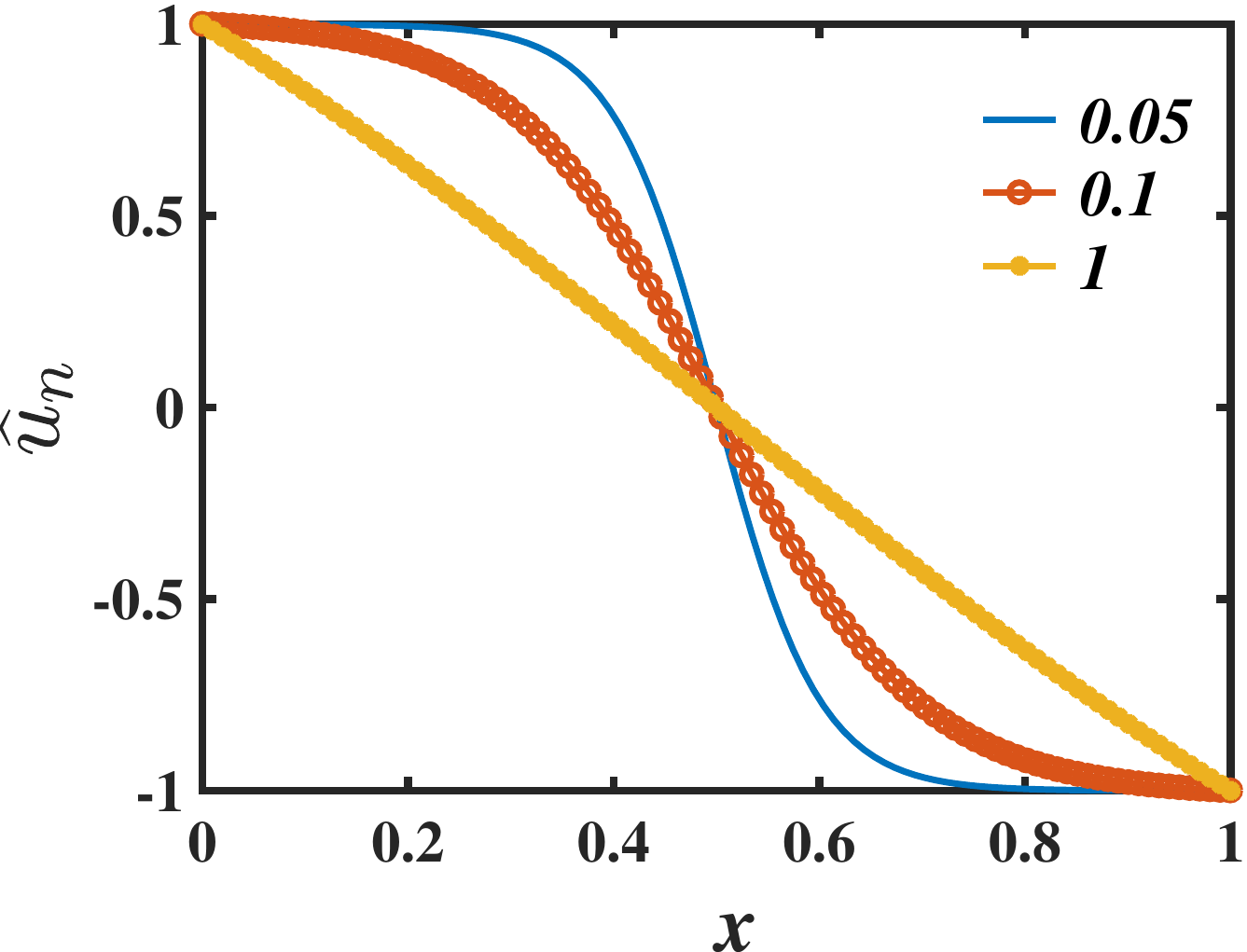}
\caption{Steady viscous Burgers' result. (Left) Histories of convergence for the error and error estimator for the traditional residual-based RBM and proposed L1-ROC. Here, $E^R$ and $E^L$ refer to the $E(n)$ in \eqref{eq:error:steadyburger} with the reduced solution $\widehat{u}_n$ constructed by following the residual-based error estimator $\Delta^R$ and L1-based importance indicator $\Delta^L$, respectively. (Middle) Distribution of selected parameters $\bmu^n$, using estimator $\Delta^R$ and $\Delta^L$, as a function of $n$. (Right) Sample RB solutions at three parameter values themselves. Note that $\Delta^L$ does not decay to zero for large $n$, but such decay is not expected or needed for this function.
}
\label{fig:steadyBurger}
\end{figure}

Error curves and the distribution of the first $N=10$ selected parameters with $\calN = 100$ are showed in Figure \ref{fig:steadyBurger}. It shows a clear exponential convergence as $n$ increases and a concentration of the selected $\mu$ values toward the lower end of the parameter domain. We note that the distributions of chosen parameters between the traditional residual-based scheme and the nascent L1-based scheme are very much similar which underscores the reliability of the new L1-ROC approach. 

\subsubsection{Nonlinear reaction diffusion equations}

 Here we consider the following cubic reaction diffusion,
\begin{equation}\label{eq:CRD}
\begin{split}
-\mu_2 \Delta u +u{(u-\mu_1)}^2 & = f(\bx) \mbox{ in } \Omega := [-1,1]\times [-1,1],\\
u & = 0 \mbox{ on } \partial \Omega.
\end{split}
\end{equation}
We take $f(\bx)=100\sin(2\pi x_1)\cos(2\pi x_2)$, and the parameter domain $\D$ is set to be $[0.2,5]\times [0.2,2]$. 
$\D$ is discretized by a $128 \times 64$ uniform tensorial grid. Denoting the step size along the $\mu_1$ direction by $h_1$, and the other by $h_2$, we specify the training set and test set as follows, 
\begin{align*}
  \Xi_{\rm train} &=  (0.2:4h_1:5) \times (0.2:4h_2:2), \\
  \Xi_{\rm test} &=   ((0.2+2h_1):4h_1:(5-2h_1)) \times ((0.2+2h_2):4h_2:(2-2h_2)),
\end{align*}
  where $(a:h:b)$ denotes an equidistant mesh over $[a,b]$ with stepsize $h$. Note in particular that the two sets defined above are disjoint. 
The nonlinear solver, based on the $5$-point stencil with $\sqrt{\calN}$ interior points at each direction of $\Omega$, for the {high fidelity truth approximation} linearizes, at the  $(\ell+1)^{\rm th}$ iteration, the equation according to
\begin{equation}
-\mu_2 \Delta u^{(\ell+1)} +g'(u^{(\ell)})u^{(\ell+1)}=g'(u^{(\ell)})u^{(\ell)}-g(u^{(\ell)},\mu_1)+f(\bx)
\label{Operator2}
\end{equation}
where $g(u;\mu_1)=u{(u-\mu_1)}^2$.

\begin{figure}[!htb]
\centering
\includegraphics[height=0.19\textheight]{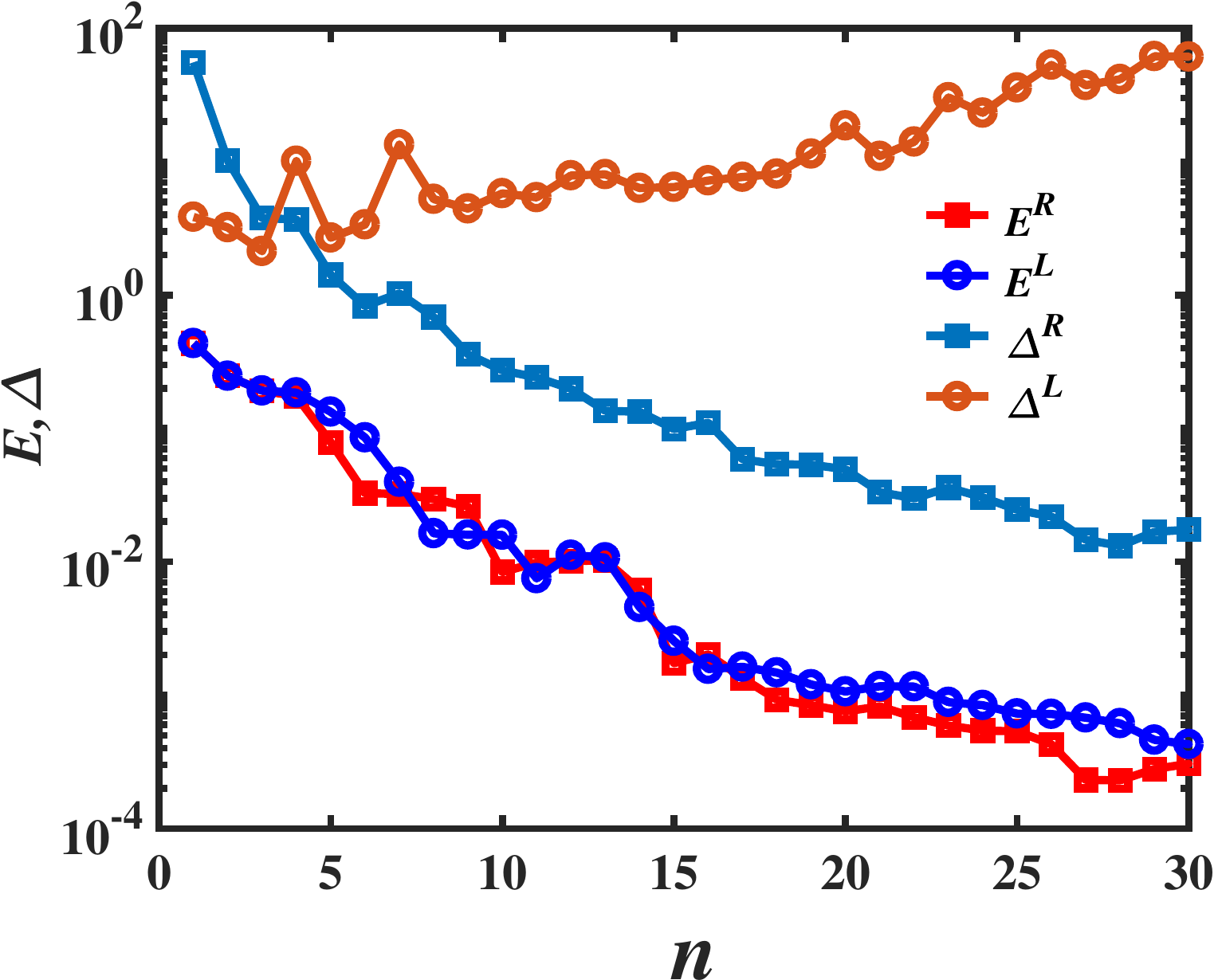}
\includegraphics[height=0.19\textheight]{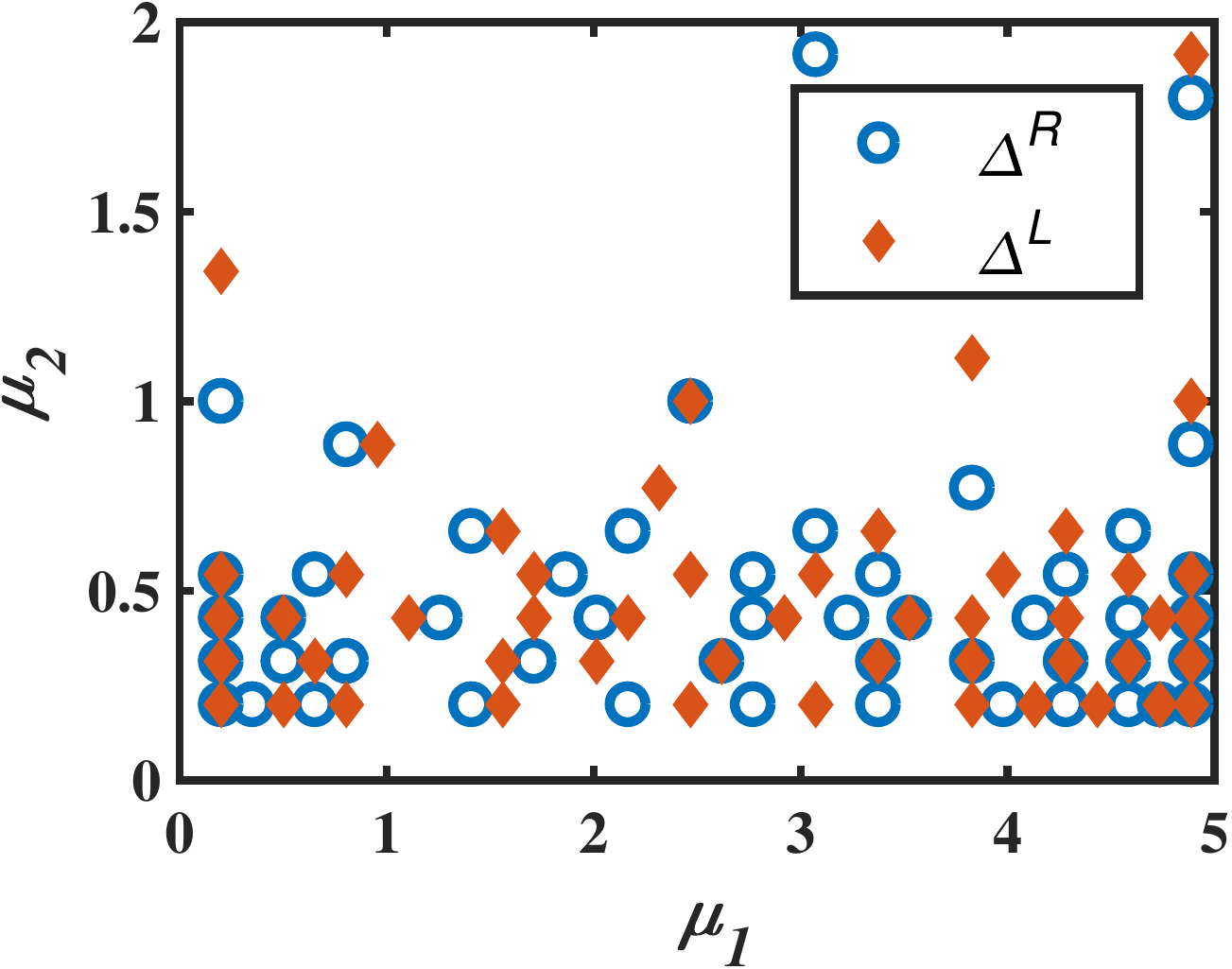}
\includegraphics[height=0.19\textheight]{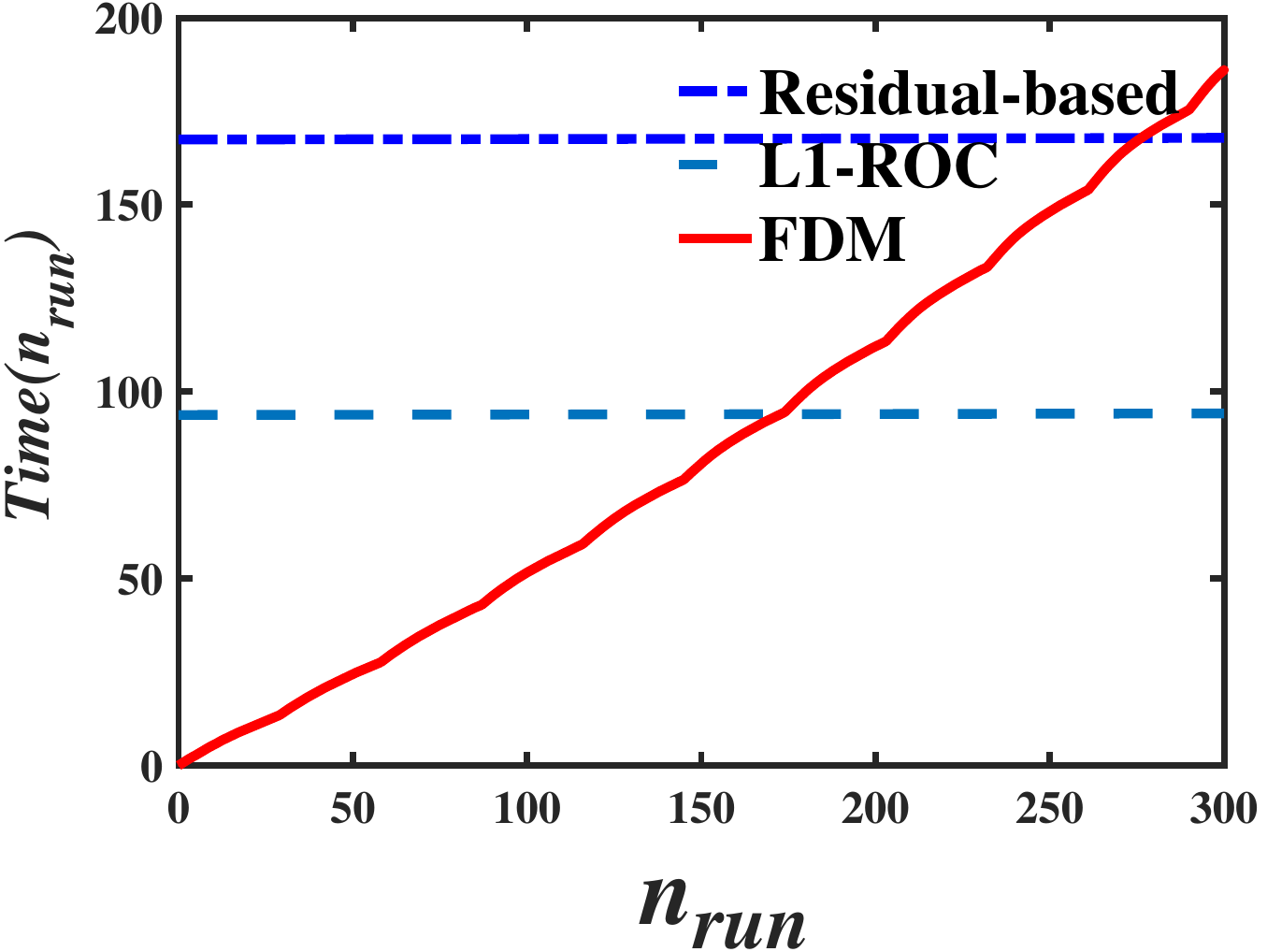}\\
 \includegraphics[width=0.45\textwidth]{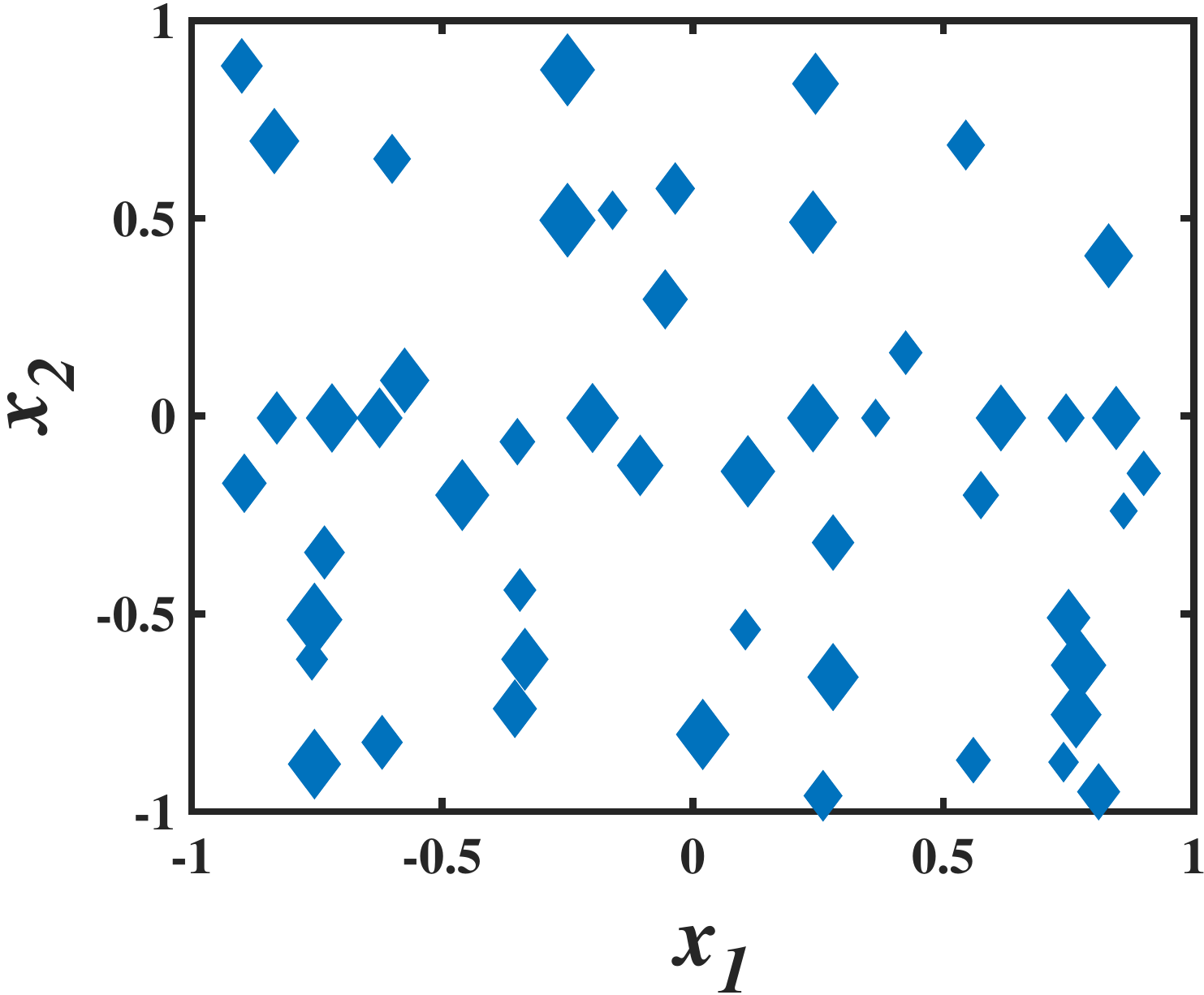}
\includegraphics[width=0.45\textwidth]{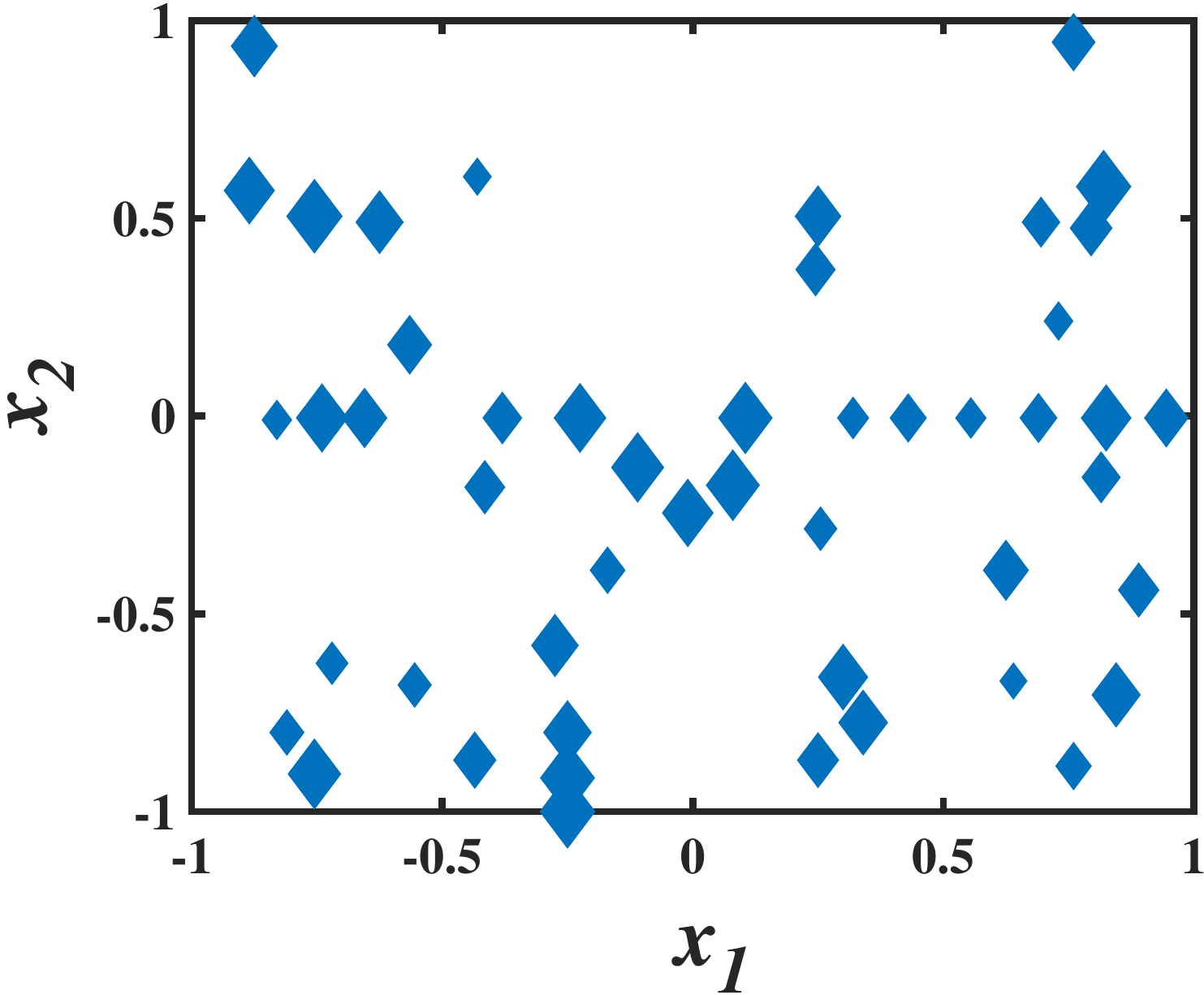}
\caption{Cubic reaction diffusion result. 
  Top row:(Left) comparison of the histories of convergence with $\sqrt{\calN} = 400$ for the errors and the error estimator for the ROC method. Here, $E^R$ and $E^L$ refer to the $E(n)$ in \eqref{eq:error:steadyburger} with the reduced solution $\widehat{u}_n$ constructed by following the residual-based error estimator $\Delta^R$ and L1-based importance indicator $\Delta^L$, respectively. (Middle) Selected $N(=40)$ parameters of the ROC method for residual-based and L1-based approaches. (Right) cumulative runtime of the FDM, the residual-based, and L1-based RBM. Bottom row: selected $40$ collocation points $X^M_s$ from solutions (Left) and $39$ collocation points $X^M_r$ from residual vector (Right).}
\label{2relativeerror}
\end{figure}

Relative errors of the RB solution $E(n)$ with $K = \sqrt{\calN} = 400$ are displayed in Figure \ref{2relativeerror} top left.  Initially, steady exponential convergence is again observed for the L1-ROC method.
The set of selected parameters are shown in Figure \ref{2relativeerror} top middle, while the collocation points are shown on the bottom row. We note again that the distributions of chosen parameters between the traditional residual-based scheme and the more nascent L1-based scheme are quite similar for this example underscoring the reliability of the L1-ROC approach.

Lastly, we showcase the vast saving of the offline time for the L1-ROC approaches. Toward that end, the comparison in cumulative computation time for the residual-based, L1-ROC, and the {high fidelity truth approximations} is shown in Figure \ref{2relativeerror} top right. 
The initial nonzero start of the L1-ROC is the amount of its offline time. 
We observe that, when $n_{\rm run}>172$, L1-ROC starts to save time in comparison to repeated runs of the truth solver. In that regard, the residual-based ROC is effective when $n_{\rm run}>276$ with $\sqrt{\calN} = 200$. 
The difference in this ``break-even'' point is because the overhead cost, devoted to calculating $\Delta_n^L$ (for L1-ROC), is significantly less than that for $\Delta_n^R$.  The latter involves (an offline-online decomposition of) the calculation of the full residual norm while the former only requires, in the L1-ROC case, obtaining an $N\times 1$ vector and evaluating its L1-norm. 
{It is worth noting that  the ``break-even'' number of runs is insensitive to $\sqrt{\calN}$}. Though L1-ROC has a much more efficient offline procedure than the  residual-based ROC, their online time for any new parameter is comparable, see Table~\ref{time2}. We observe that the L1-ROC method accelerates the iterative truth solver by $2000 \sim 50000$ times. The results also confirm that time consumption of the online ROC methods is independent of $K=\sqrt{\calN}$. 
In order to demonstrate the time savings more intuitively, we present the online calculation time for the different algorithms in two different parameter regimes. The first regime is when $\mu_1$ is large and $\mu_2$ small, in particular we choose $\mu_1=4.55, \mu_2=0.42$. The second regime has the relative sizes reversed. The reduced solver requires $27$ iterations for the nonlinear system in the first regime, while only requiring $8$ iterations in the second regime. Therefore, the full-order time consumption seems very different. However, Table \ref{time2} does indicate a speedup range of $3000 \sim 17000$ when $\sqrt{\calN}=400,\, 800$.

\begin{table}[!htb]
	\centering
		\begin{tabular}{ccccc}
		\hline
$(\mu_1, \mu_2)$ &~~$K$~~&Residual-based ROC &  ~ L1-ROC~ &~Direct FDM~~~~~  \\ \hline
\multirow{3}{*}{$(4.55, 0.42)$} &200	&0.003150  &0.003159 &2.310034 \\ 
&400	&  0.003067  &0.003136 & 11.779558 \\ 
&800	&  0.003258  &0.003162 &53.727031 \\ \hline
\multirow{3}{*}{$(1, 1.82)$} &200	&0.001125  &0.001060 &0.662095 \\ 
&400	& 0.001141 & 0.001205 &3.338956  \\ 
&800	&0.001207  &0.001261& 15.173460 \\ \hline
	\end{tabular}
        \caption{Online computational times (seconds) with different grid sizes $K$, when $N=40$. }
	\label{time2}
\end{table}

\begin{remark}
The L1-ROC works well for other nonlinear convection diffusion reaction equations. For example, we tested the dimensionless nonlinear nonaffine Poisson-Boltzmann equation 
\begin{align}\label{eq:PB}
D\nabla^2 u  = \sinh u+ g(\bx), \mbox{ with } g(\bx)= \exp[-50({(x_1-0.2)}^2+{(x_2+0.1)}^2)] 
\end{align}
modeling a source distribution centered at $(0.2,-0.1)$. The parameters are diffusion coefficient $D$ and the voltage differential $V$ at the boundary. The authors have previously designed a RBM for this equation \cite{JCX2018}. However, due to the desire to avoid applying EIM directly, we observed limited speedup (less than one order of magnitude). With L1-ROC, we achieved a speedup factor of up to four orders of magnitude, see Table \ref{table:pb}.
This significant progress underscores the power of the L1-ROC approach. In addition, we tested an equation with nonlinear convection term 
\begin{equation}
-\mu_2 \Delta u + u \left( \lVert \nabla u \rVert + \mu_1 \right)^{1.5} = f(x).
\label{eq:conv}
\end{equation}
The L1-ROC works equally effective as well, see Table \ref{table:convection}.

\begin{table}[!htb]
	\centering
%		\begin{tabular}{lrrr}
			\begin{tabular}{cccc}
		\hline
~~$K$~~&Residual-based ROC &  ~ L1-ROC~&~Direct FDM~~~~~  \\ \hline
200	&0.000678 &0.000688&   1.439812 \\ 
400	&0.000770  &0.000646 &  6.492029  \\ 
800	&0.000728    &0.000625&  33.722112 \\ \hline
	\end{tabular}
	\caption{Online computational times (in seconds) for the Poisson-Boltzmann equation \eqref{eq:PB} at different grid sizes $K$, when $V = 3.85, D= 0.15^2$, $N=30$.}
	\label{table:pb}
\end{table}

\begin{table}[!htb]
	\centering
				\begin{tabular}{cccc}
		\hline
~~$K$~~&Residual-based ROC &  ~ L1-ROC~ & ~Direct FDM~~~~  \\ \hline
200	&0.000422 &0.000428 &0.569732 \\ 
400	& 0.000397 & 0.000410& 2.838783  \\ 
800	&0.000424  &0.000425&12.582593 \\ \hline
	\end{tabular}
	\caption{Online computational times (seconds) for the nonlinear convection diffusion equation \eqref{eq:conv} at different grid sizes $K$, when $N=20, \mu_1 = 32, \mu_2=3$.}
	\label{table:convection}
\end{table}

\end{remark}

%\clearpage

\subsubsection{Numerical comparison with POD and random generation} 

\begin{figure}[!htb]
\centering
\includegraphics[width=0.32\textwidth]{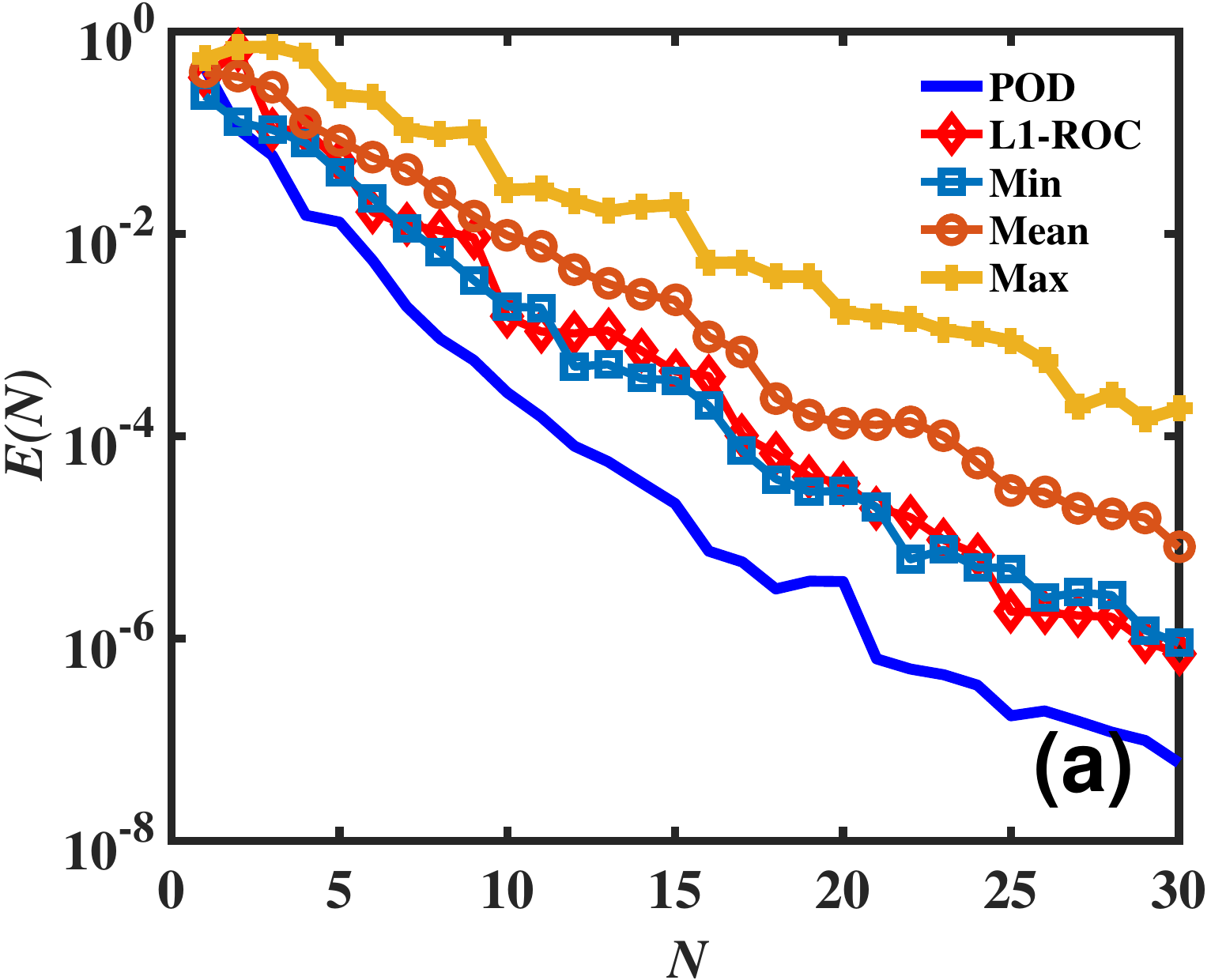}
\includegraphics[width=0.32\textwidth]{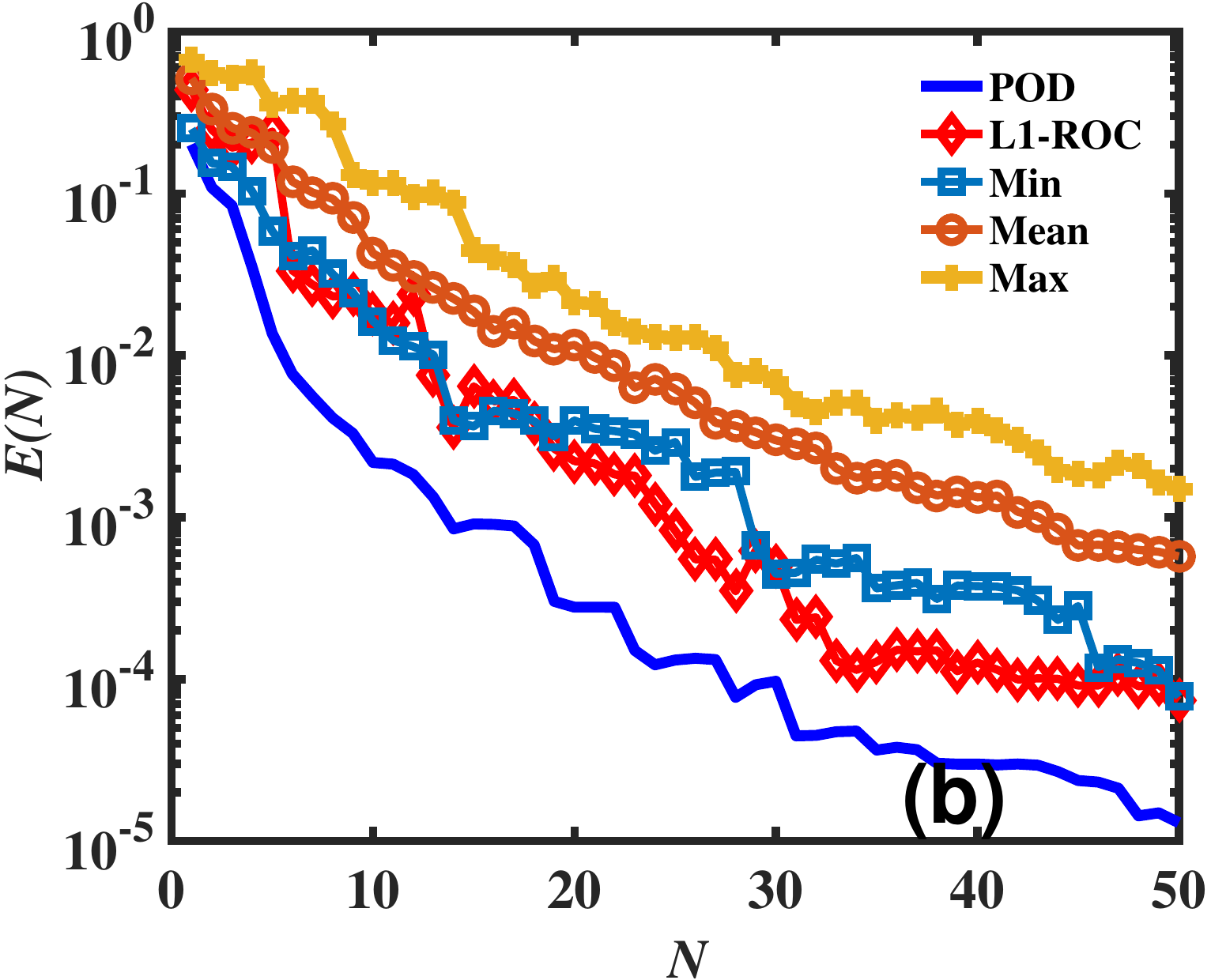}
\includegraphics[width=0.34\textwidth]{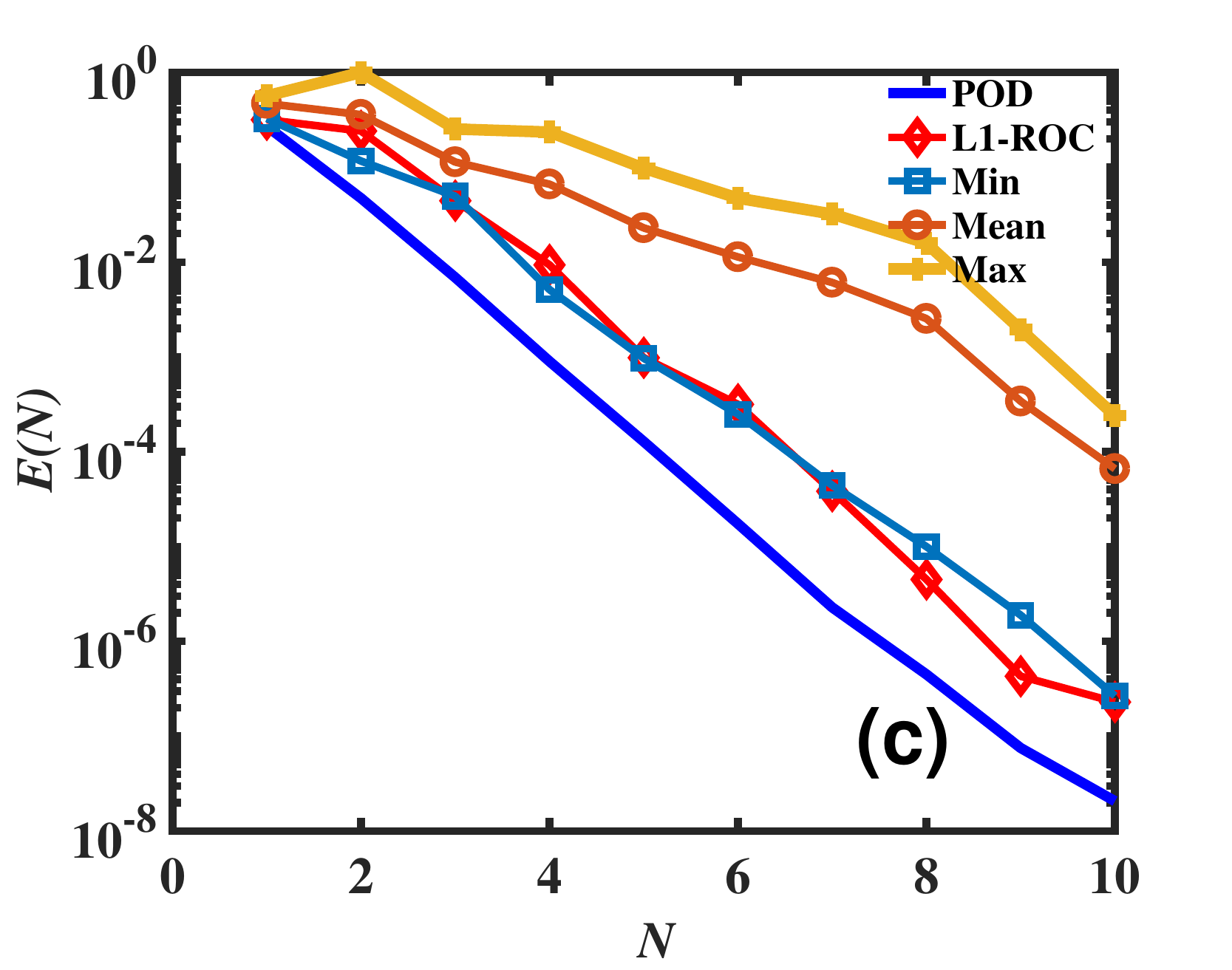}
  \caption{Convergence comparison for the L1-ROC, exhaustive POD and (best, median, and worst cases of) random generation approaches. (a) Poisson-Boltzmann equation \eqref{eq:PB} with $\sqrt{\calN}=400$, (b) cubic reaction diffusion \eqref{eq:CRD} with $\sqrt{\calN}=400$, (c) steady viscous Burgers' equation \eqref{eq:burgers} with $\sqrt{\calN}=100$.}
\label{1compareerror}
\end{figure}

To further establish numerically the reliability of the L1-ROC algorithm, we compare it with two alternative methods of building the reduced basis space. On one end, the proper orthogonal decomposition (POD) \cite{BerkoozHolmesLumley1993, Kunisch_Volkwein_POD, WillcoxPeraire2002, LiangPOD} based on an exhaustive selection of snapshots (i.e. we include all solutions $u^\N(\bmu)$ for $\bmu \in \Xi_{\rm train}$)  produces the best reduced solution space and thus the most accurate, albeit costly, surrogate solution. We note that this version of POD only serves as reference and is in general not feasible as the full solution ensemble must be generated.  On the other end, a random selection of $N$ parameters as our RB snapshots is a fast but crude method. 
Comparison results of three steady-state test problems are shown in Figure \ref{1compareerror} with FDM points per dimension $\sqrt{\calN}$ set to be $400$  for first two cases (results with different $\sqrt{\calN}$ are similar) and $\sqrt{\calN}=100$ for the third case. Not surprisingly, the exhaustive POD is the most accurate. 
Our L1-ROC is one order of magnitude worse than POD, but in fact slightly better or comparable to the the best possible random generation. It is roughly one order of magnitude better than the median performance of random generations.

\subsection{Time dependent nonlinear problems}
\label{numerics:timedep}
In this section, we test the time-dependent equations corresponding to stationary problems in the last section, namely viscous Burgers' and cubic reaction diffusion equations.

\subsubsection{Viscous Burgers' equation}

\begin{figure}[!htb]
\centering
\includegraphics[width=0.47\textwidth]{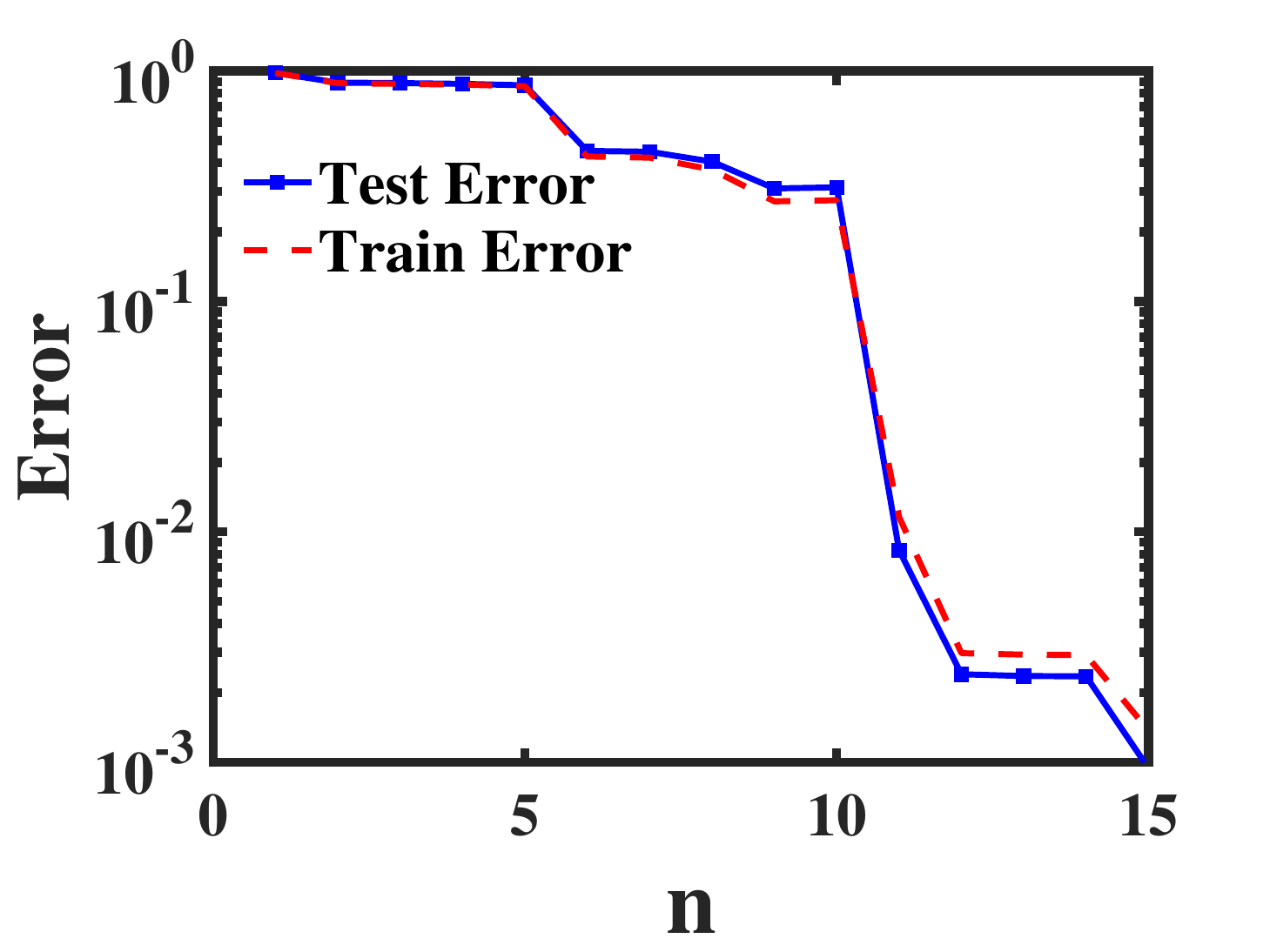}
\includegraphics[width=0.45\textwidth]{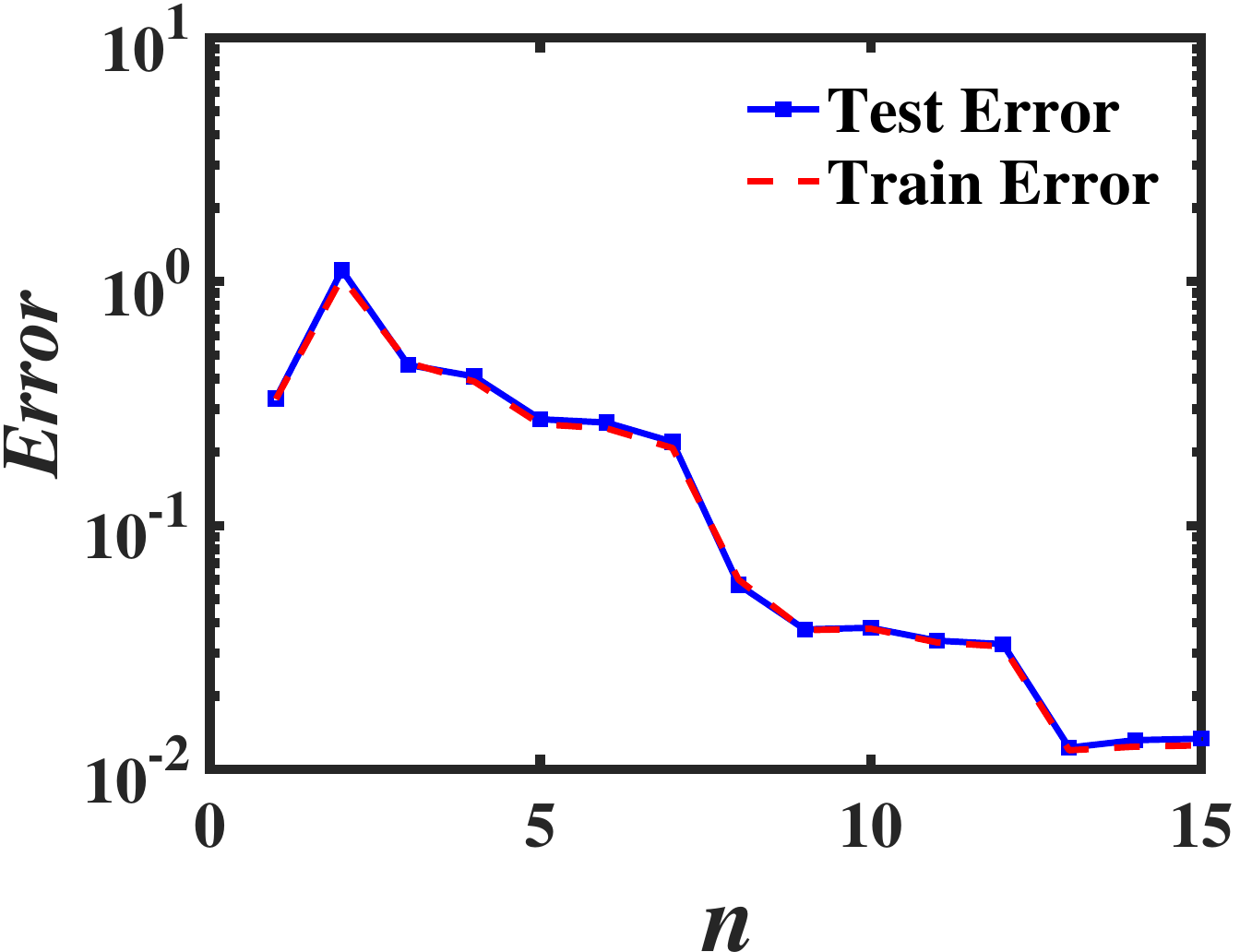}\\
\includegraphics[height=0.18\textheight]{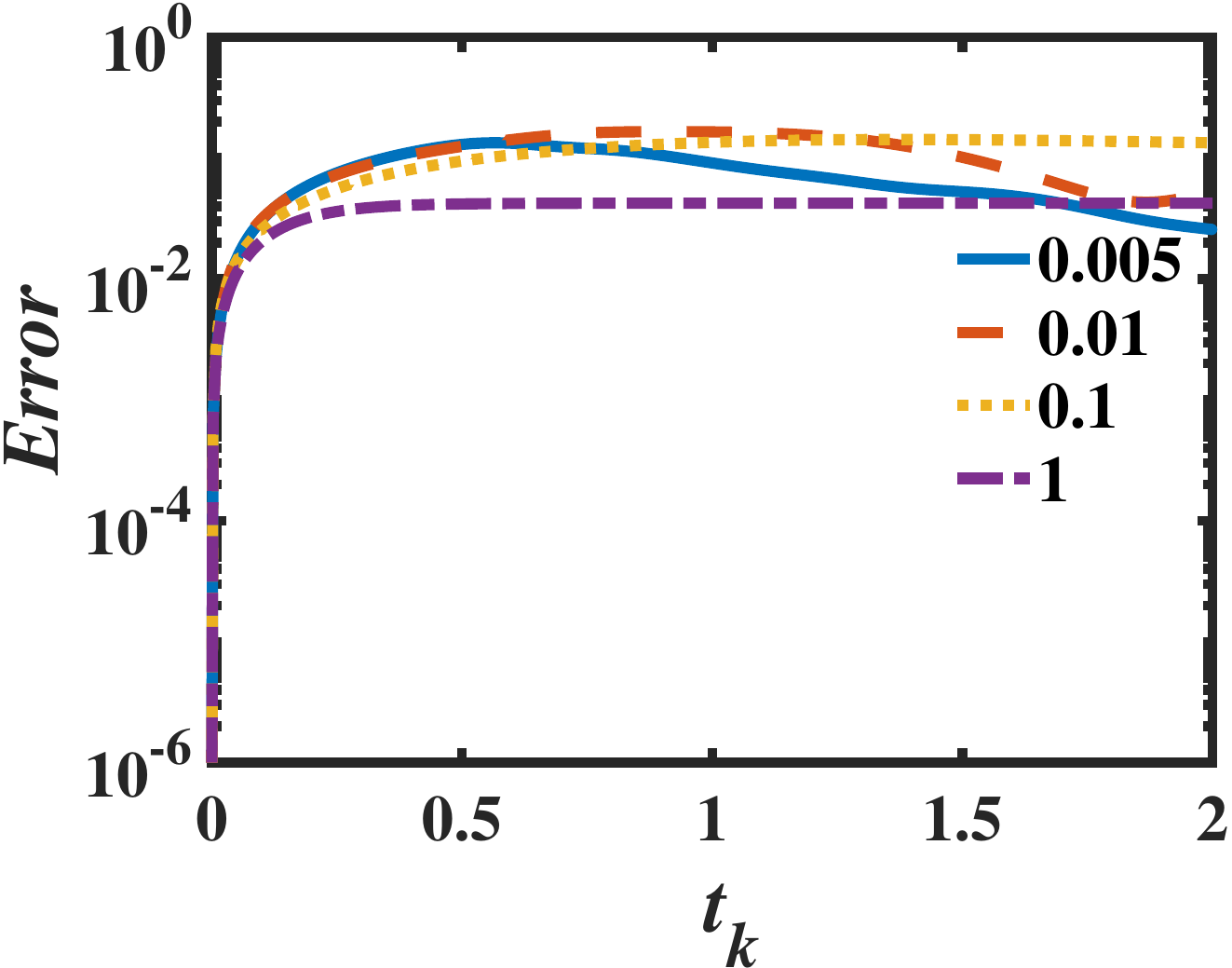}
\includegraphics[height=0.19\textheight]{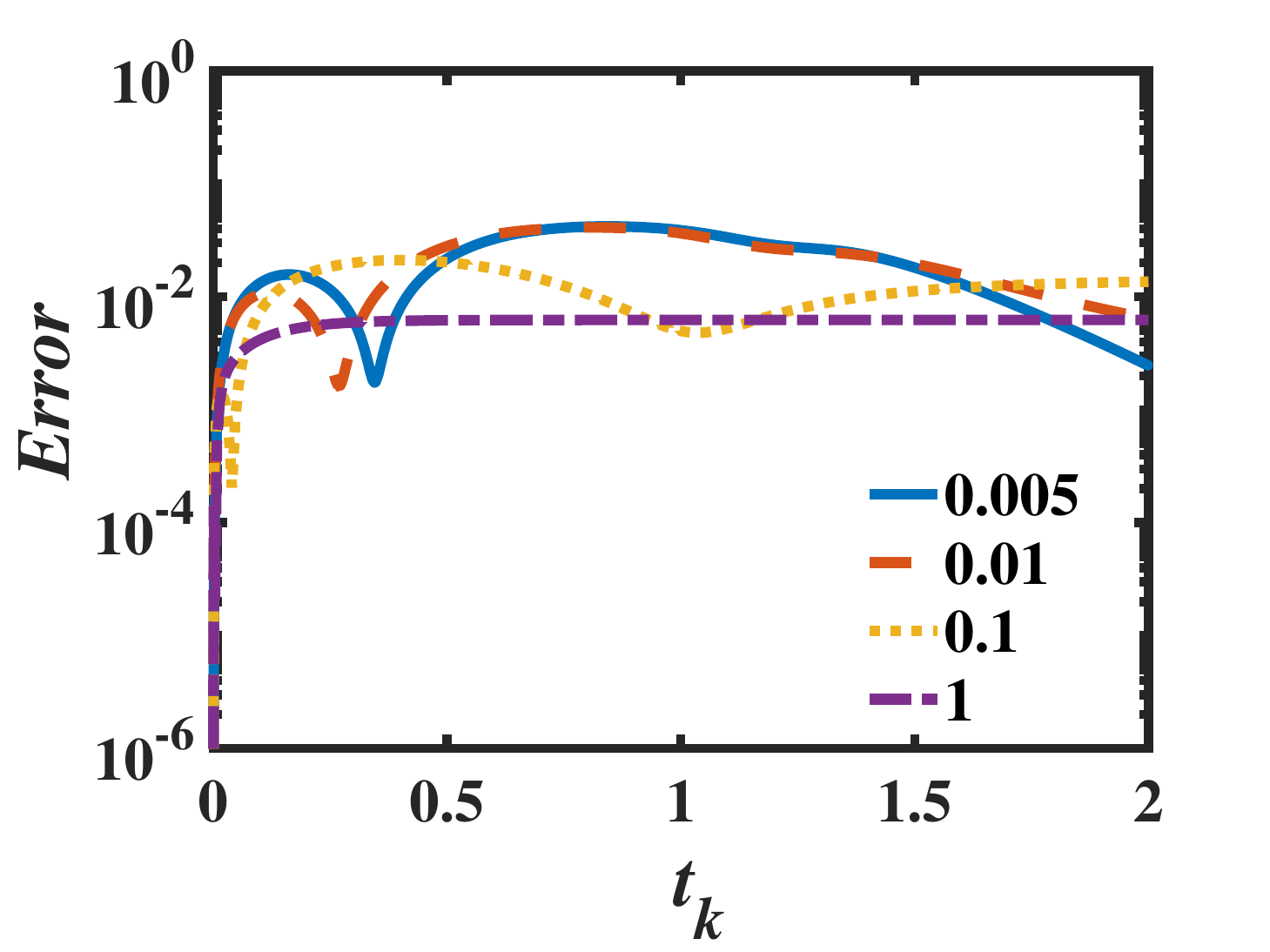}
\includegraphics[height=0.19\textheight]{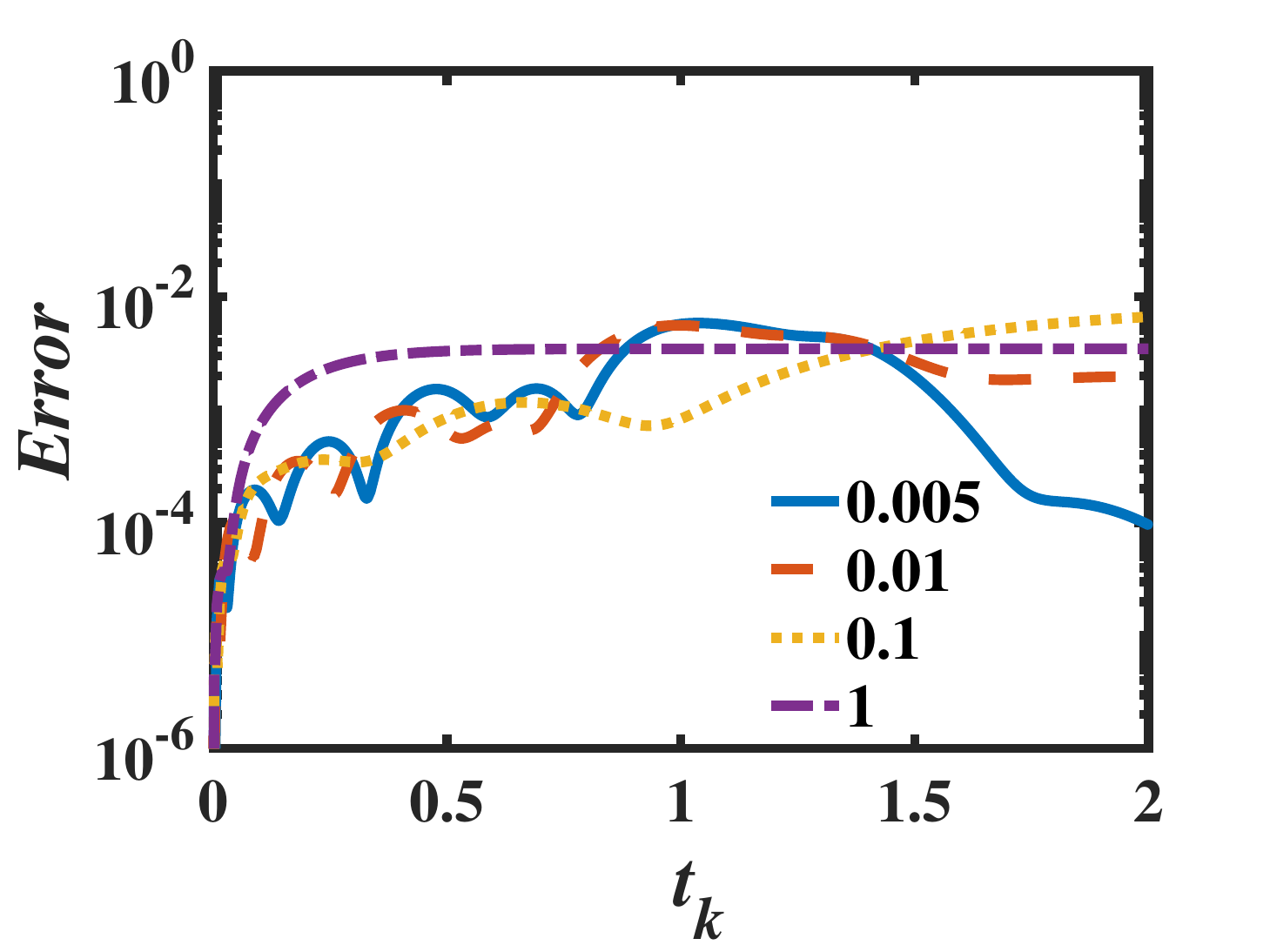}
\caption{Transient viscous Burgers' result. On the top row are the error curves of L1-ROC with $N=15$ basis elements for the setup in \cite{peherstorfer2019sampling} (left) and \cite{nguyen2009reduced} (right). Plotted at the bottom are the actual $L^2$ error, $||u^\N(:,t_k;\bmu) -u_N(:,t_k;\bmu)||$ as a function of discrete time $t_k$. The left, center and right plots show $N=5, 10, 15$, respectively, each for parameter values $\bmu = 0.005,0.01, 0.1,1$ with the setup as in \cite{nguyen2009reduced}. 
}
\label{fig:ben}
\end{figure}

We test the viscous Burgers' equation adopting settings similar to \cite{peherstorfer2019sampling, nguyen2009reduced}
\begin{equation}
\begin{split}
u_t + u u_x & = \mu u_{xx} + f(x),~ (x, t, \mu) \in (0,1) \times (0,1] \times \calD,\\
u(x,t=0; \mu) & =0,\\
u(0, t; \mu) = \alpha, \,\, & \,\, u(1, t; \mu) = \beta.
\end{split} 
\end{equation} 
The authors of \cite{peherstorfer2019sampling} takes $\calD = [0.1,1], f = 0, T = 1, \Delta t =10^{-4}, (\alpha, \beta) = (-1, 1)$ and monitor the average error in a Frobenius norm-based metric,
\begin{align*}
  \textrm{Error} = \frac{1}{m_{test}}\sum_{i=1}^{m_{test}} \frac{||u(\cdot,\cdot;\bmu) -\widehat{u}(\cdot,\cdot;\bmu)||_F}{||u(\cdot,\cdot;\bmu)||_F},~~ 
  \| v(\cdot,\cdot) \|^2_F &\coloneqq \sum_{\bx \in X^\N, t_i \in {\mathcal T}_f} v(\bx,t_i)^2
\end{align*}
while the authors in \cite{nguyen2009reduced} set $\calD = [0.005,1], f = 1, T = 2, \Delta t = 2 \cdot 10^{-6}, (\alpha, \beta) = (0,0)$ and observe the error in $L^2$. We investigate L1-ROC results from both of these setups. The results are showed in Figure \ref{fig:ben}. These results are similar to those of \cite{peherstorfer2019sampling, nguyen2009reduced}. However, we note that they come at a much smaller computational expense.

\subsubsection{Nonlinear reaction diffusion problems}

Next, we consider accordingly the following time dependent nonlinear reaction diffusion equation,
\begin{equation}
\begin{split}
u_t -\mu_2 \Delta u +u{(u-\mu_1)}^2 & = f(\bx), \mbox{ in } \Omega=[-1,1]\times [-1,1],\\
u & = 0 \mbox{ on } \partial \Omega, \\
u(\bx, t = 0) & = u_0(\bx).
\end{split}
\end{equation}
Here $f(\bx)=100\sin(2\pi x_1)\cos(2\pi x_2)$, and $[\mu_1,\mu_2] \in \calD := [1,5]\times [0.2,1]$. 
The parameter space $\calD$ is discretized by a $128 \times 32$ uniform tensorial grid. Denoting the step size along the $\mu_1$ direction by $h_1$, and the other by $h_2$, we specify the training and test sets as follows, 
\begin{align*}
  \Xi_{\rm train} &=  (1:8h_1:5) \times (0.2:2h_2:1), \\
  \Xi_{\rm test} &=   ((1+2h_1):4h_1:(5-2h_1)) \times ((0.2+h_2):4h_2:(1-h_2)),
\end{align*}
For the truth approximation, we use backward Euler for time marching and the same nonlinear spatial solver as the steady-state case \eqref{Operator2}.

\begin{figure}[!htb]
\centering
\includegraphics[width=0.45\textwidth]{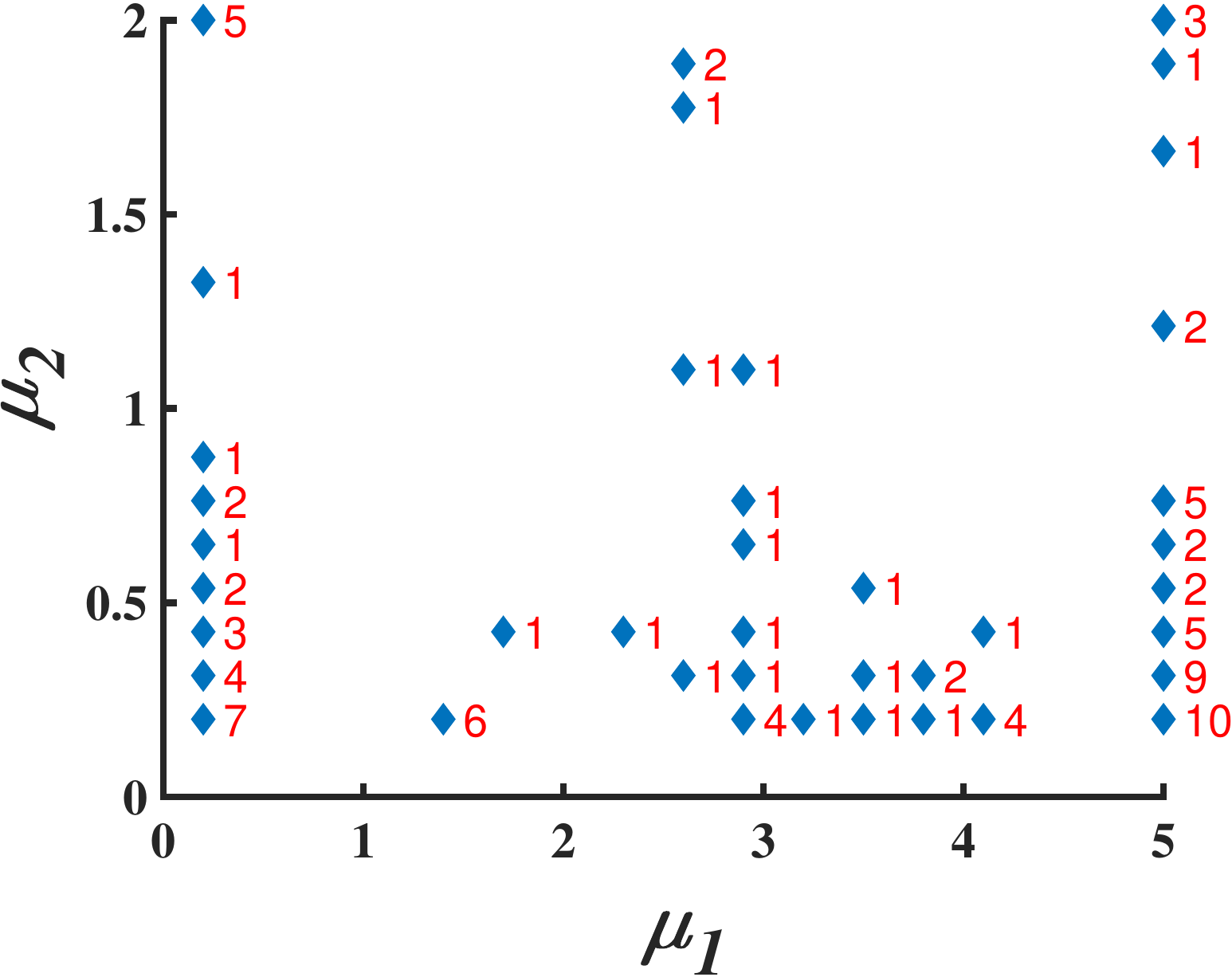}
\includegraphics[width=0.54\textwidth]{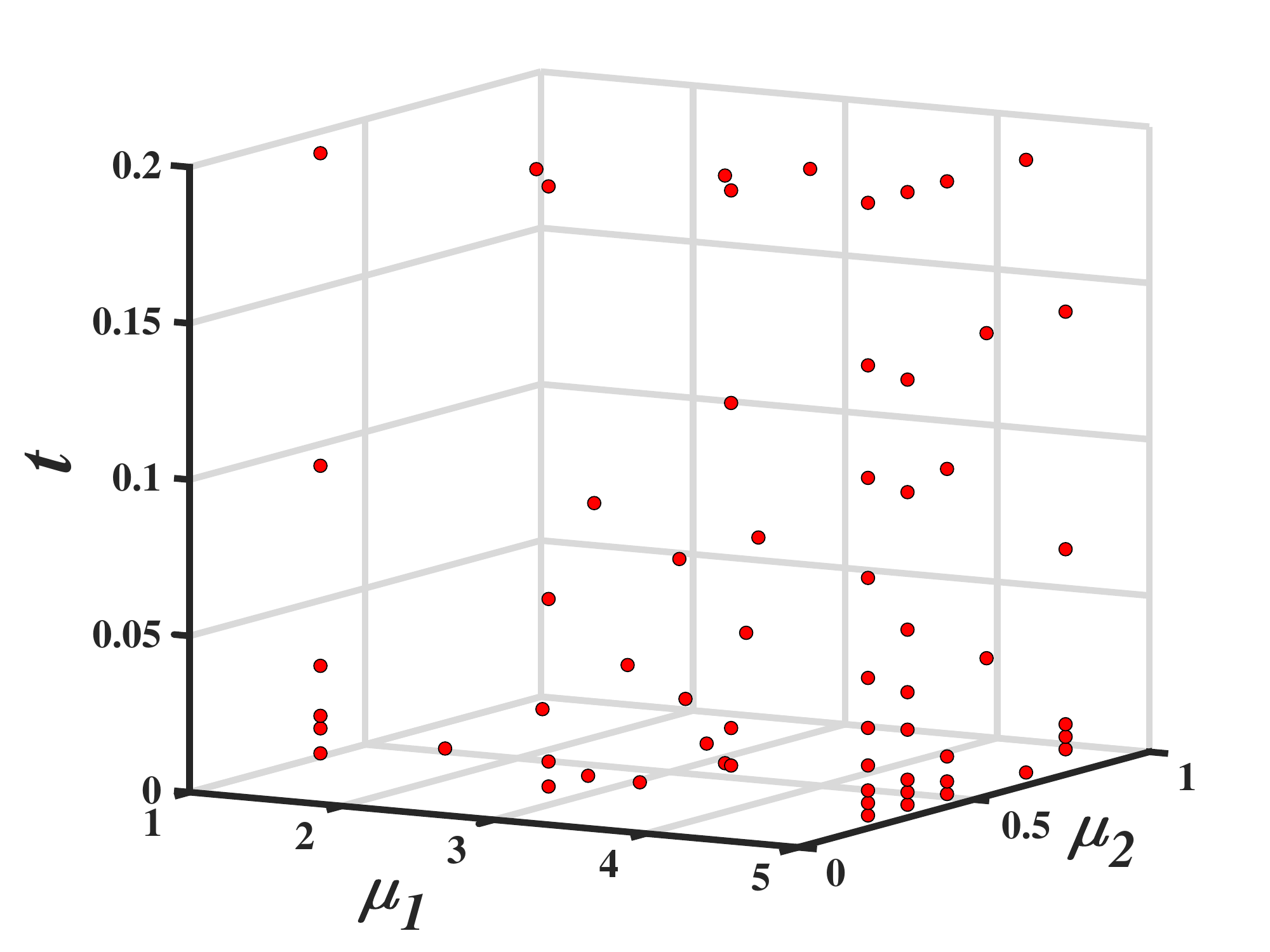}
\includegraphics[height=0.18\textheight]{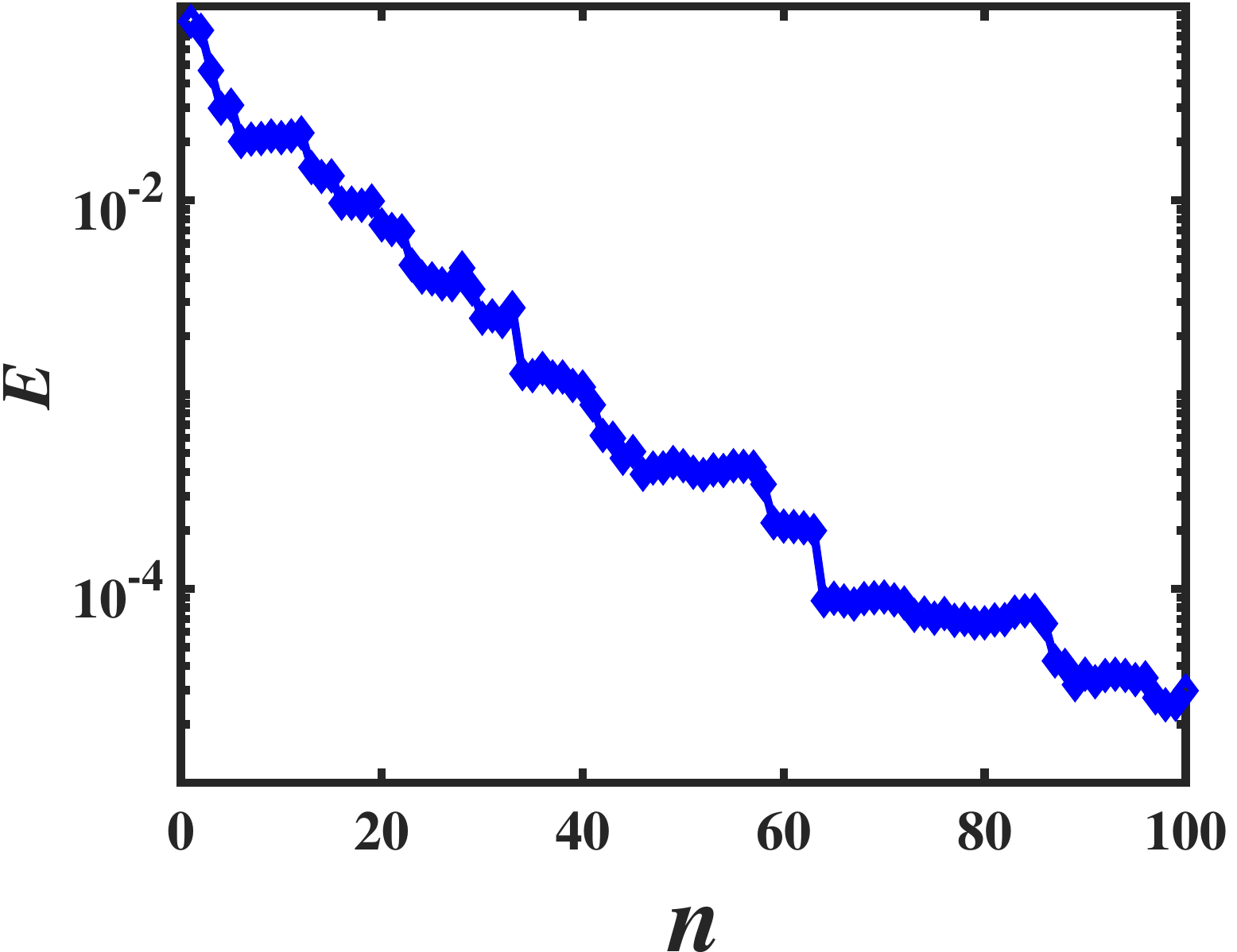}
\includegraphics[height=0.19\textheight]{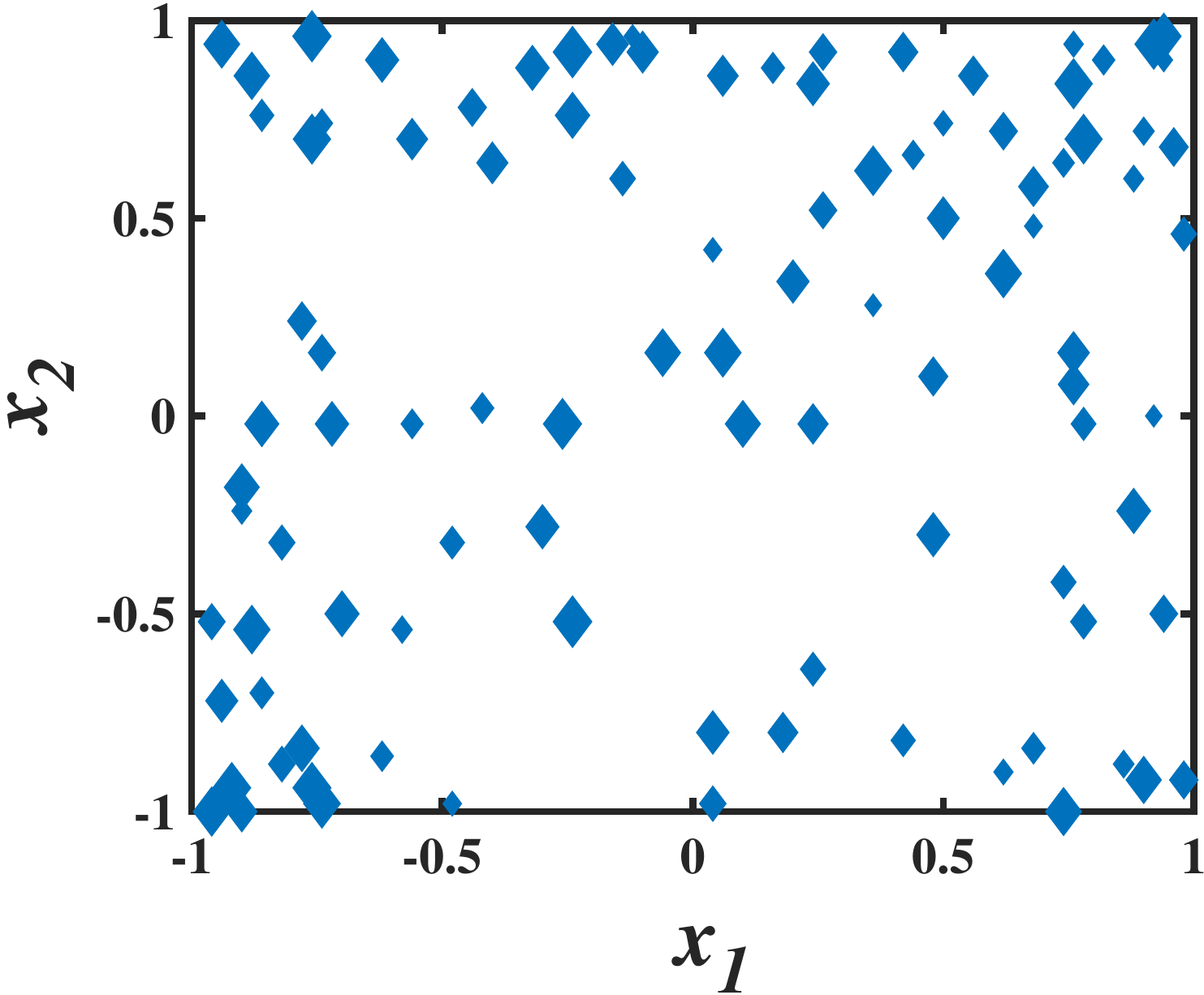}
\includegraphics[height=0.19\textheight]{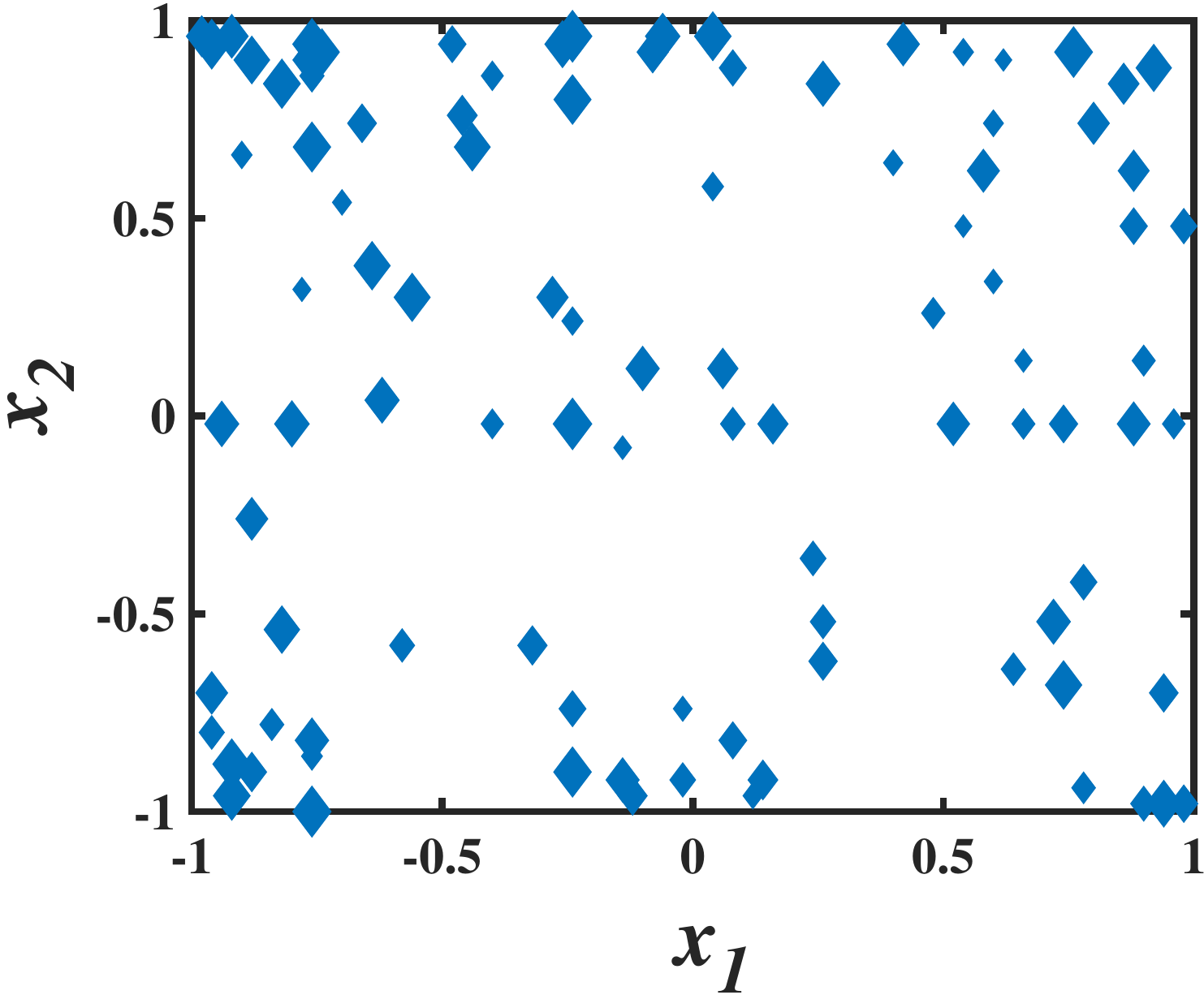}
\caption{Transient cubic reaction diffusion result. Top Left: Selected parameters when $N_{max}=100$. The number means corresponding parameter is selected at many different time nodes. Top Right: a three-dimensional view of the selected parameters. Error curves of L1-ROC algorithm, and collocation points from solutions and residuals are shown at the bottom row from left fo right respectively. 
}
\label{fig:timecubic}
\end{figure}

We report the $\bmu$-component of  the parameter values selected by L1-ROC in Figure \ref{fig:timecubic} (top). Note that the RB space is built from the snapshots
\[
\left\{u(t^1_{\bmu^n}, \cdot; \bmu^n), \dots, u(t_{\bmu^n}^{k_{\bmu^n}}, \cdot; \bmu^n)\right\}_{n=1}^N.
\]
That is, for each distinct parameter value $\bmu^n$ chosen by L1-ROC, there are $k_{\bmu^n} \ge 1$ time level  snapshots $\{t_{\bmu^n}^1, \dots, t_{\bmu^n}^{k_{\bmu^n}}\} \subset \{t_0, t_1, \dots, t_{\calN_t}\}$. The red number by each $\bmu$ values in the left pane denotes this $k_{\bmu^n}$. It is interesting to note that, consistent with the tendency of RBM selecting boundary values of the parameter domain, our L1-ROC tends to select multiple snapshots along time for the selected parameters along the boundary of the parameter domain. 
The right pane is the corresponding 3D-image of the left. 
The bottom row of Figure \ref{fig:timecubic} shows the L1-ROC error curve, which shows clear exponential convergence, and collocation points in the physical domain.

\section{Conclusion}
\label{sec:conclusion}

This paper proposes a novel reduced over-collocation method, dubbed L1-ROC, for efficiently solving parametrized nonlinear and nonaffine PDEs. By integrating EIM technique on the solution snapshots and well-chosen residuals, the collocation philosophy, and the simplicity of the L1-based importance indicator that is extended to time-dependent problems, L1-ROC has online computational complexity independent of the degrees of freedom of the underlying FDM, and  furthermore immune from the number of EIM expansion terms. 
This expansion would have otherwise significantly degraded the efficiency of a traditional RBM when applied to the nonaffine and nonliner terms in the equation. The lack of such precomputations of nonlinear and nonaffine terms makes the method dramatically faster offline and online, and significantly simpler to implement than any existing RBM. For future directions, we plan to apply L1-ROC to systems of equations resulting from CFD systems with more complicated nonlinear and nonaffine terms. A deeper understanding on the theory of this L1-ROC algorithm is also in our consideration.

\bibliographystyle{abbrv}
\bibliography{rbmbib}

\begin{thebibliography}{10}

\bibitem{Almroth1978}
B.~O. Almroth, P.~Stern, and F.~A. Brogan.
\newblock Automatic choice of global shape functions in structural analysis.
\newblock {\em AIAA J.}, 16(5):525--528, 1978.

\bibitem{astrid2004fast}
P.~Astrid.
\newblock Fast reduced order modeling technique for large scale ltv systems.
\newblock In {\em Proceedings of the 2004 American control conference},
  volume~1, pages 762--767. IEEE, 2004.

\bibitem{Astrid2008MissingPoint}
P.~{Astrid}, S.~{Weiland}, K.~{Willcox}, and T.~{Backx}.
\newblock Missing point estimation in models described by proper orthogonal
  decomposition.
\newblock {\em IEEE Transactions on Automatic Control}, 53(10):2237--2251,
  2008.

\bibitem{Barrault2004}
M.~Barrault, Y.~Maday, N.~C. Nguyen, and A.~T. Patera.
\newblock An 'empirical interpolation' method: {A}pplication to efficient
  reduced-basis discretization of partial differential equations.
\newblock {\em C. R. Math.}, 339(9):667--672, 2004.

\bibitem{BenaceurEhrlacherErnMeunier2018}
A.~Benaceur, V.~Ehrlacher, A.~Ern, and S.~Meunier.
\newblock A progressive reduced basis/empirical interpolation method for
  nonlinear parabolic problems.
\newblock {\em SIAM J. Sci. Comput.}, 40(5):A2930--A2955, 2018.

\bibitem{BennerGugercinWillcox2015}
P.~Benner, S.~Gugercin, and K.~Willcox.
\newblock A survey of projection-based model reduction methods for parametric
  dynamical systems.
\newblock {\em SIAM review}, 57(4):483--531, 2015.

\bibitem{bos2004accelerating}
R.~Bos, X.~Bombois, and P.~Van~den Hof.
\newblock Accelerating large-scale non-linear models for monitoring and control
  using spatial and temporal correlations.
\newblock In {\em Proceedings of the 2004 American Control Conference},
  volume~4, pages 3705--3710. IEEE, 2004.

\bibitem{CarlbergBaroneAntil2017}
K.~Carlberg, M.~Barone, and H.~Antil.
\newblock Galerkin v. least-squares {P}etrov-{G}alerkin projection in nonlinear
  model reduction.
\newblock {\em J. Comput. Phys.}, 330:693--734, 2017.

\bibitem{CarlbergBouMoslehFarhat2011}
K.~Carlberg, C.~Bou-Mosleh, and C.~Farhat.
\newblock Efficient non-linear model reduction via a least-squares
  {P}etrov-{G}alerkin projection and compressive tensor approximations.
\newblock {\em Int J Numer Methods Eng.}, 86(2):155--181, 2011.

\bibitem{CARLBERG2013623}
K.~Carlberg, C.~Farhat, J.~Cortial, and D.~Amsallem.
\newblock The {GNAT} method for nonlinear model reduction: Effective
  implementation and application to computational fluid dynamics and turbulent
  flows.
\newblock {\em J. Comput. Phys.}, 242:623 -- 647, 2013.

\bibitem{Casenave2014_M2AN}
F.~Casenave, A.~Ern, and T.~Leli\`evre.
\newblock Accurate and online-efficient evaluation of the {\it a posteriori}
  error bound in the reduced basis method.
\newblock {\em ESAIM Math. Model. Numer. Anal.}, 48(1):207--229, 2014.

\bibitem{ChaturantabutSorensen2010}
S.~Chaturantabut and D.~C. Sorensen.
\newblock {Nonlinear model reduction via discrete empirical interpolation}.
\newblock {\em SIAM J. Sci. Comput.}, 32(5):2737--2764, 2010.

\bibitem{ChenGottlieb2013}
Y.~Chen and S.~Gottlieb.
\newblock Reduced collocation methods: {Reduced basis methods in the
  collocation framework}.
\newblock {\em J. Sci. Comput.}, 55(3):718--737, 2013.

\bibitem{ChenGottliebMaday}
Y.~Chen, S.~Gottlieb, and Y.~Maday.
\newblock {Parametric analytical preconditioning and its applications to the
  reduced collocation methods.}
\newblock {\em C. R. Acad. Sci. Paris, Ser. I}, 352:661--666, 2014.

\bibitem{chen2010certified}
Y.~Chen, J.~S. Hesthaven, Y.~Maday, and J.~Rodriguez.
\newblock Certified reduced basis methods and output bounds for the harmonic
  {M}axwell's equations.
\newblock {\em SIAM J. Sci. Comput.}, 32(2):970--996, 2010.

\bibitem{chen2012certified}
Y.~Chen, J.~S. Hesthaven, Y.~Maday, J.~Rodriguez, and X.~Zhu.
\newblock Certified reduced basis method for electromagnetic scattering and
  radar cross section estimation.
\newblock {\em Comput. Methods Appl. Mech. Engrg.}, 233-236:92--108, 2012.

\bibitem{JiangChenNarayan2019}
Y.~Chen, J.~Jiang, and A.~Narayan.
\newblock {A robust error estimator and a residual-free error indicator for
  reduced basis methods}.
\newblock {\em Computers \& Mathematics with Applications}, 77:1963--1979,
  2019.

\bibitem{deparis2009reduced}
S.~Deparis and G.~Rozza.
\newblock Reduced basis method for multi-parameter-dependent steady
  {Navier--Stokes} equations: applications to natural convection in a cavity.
\newblock {\em J. Comput. Phys.}, 228(12):4359--4378, 2009.

\bibitem{Dimitriu2017}
G.~Dimitriu, R.~{\c{S}}tef{\u{a}}nescu, and I.~M. Navon.
\newblock Comparative numerical analysis using reduced-order modeling
  strategies for nonlinear large-scale systems.
\newblock {\em Journal of Computational and Applied Mathematics}, 310:32--43,
  2017.

\bibitem{Everson1995}
R.~Everson and L.~Sirovich.
\newblock Karhunen--loeve procedure for gappy data.
\newblock {\em JOSA A}, 12(8):1657--1664, 1995.

\bibitem{Fritzen}
F.~Fritzen, B.~Haasdonk, D.~Ryckelynck, and S.~Sch{\"o}ps.
\newblock An algorithmic comparison of the hyper-reduction and the discrete
  empirical interpolation method for a nonlinear thermal problem.
\newblock {\em Mathematical and Computational Applications}, 23(1):8, 2018.

\bibitem{BerkoozHolmesLumley1993}
P.~H. G.~Berkooz and J.~Lumley.
\newblock The proper orthogonal decomposition in the analysis of turbulent
  flows.
\newblock {\em Ann. Rev. Fluid Mech.}, 25(1):539--575, 1993.

\bibitem{galbally2010non}
D.~Galbally, K.~Fidkowski, K.~Willcox, and O.~Ghattas.
\newblock Non-linear model reduction for uncertainty quantification in
  large-scale inverse problems.
\newblock {\em Int J Numer Methods Eng.}, 81(12):1581--1608, 2010.

\bibitem{grepl2012}
M.~A. Grepl.
\newblock Certified reduced basis methods for nonaffine linear time-varying and
  nonlinear parabolic partial differential equations.
\newblock {\em Mathematical Models and Methods in Applied Sciences},
  22(03):1150015, 2012.

\bibitem{grepl2007efficient}
M.~A. Grepl, Y.~Maday, N.~Nguyen, and A.~T. Patera.
\newblock Efficient reduced-basis treatment of nonaffine and nonlinear partial
  differential equations.
\newblock {\em ESAIM-Math. Model. Numer. Anal.}, 41(3):575--605, 2007.

\bibitem{grepl2005posteriori}
M.~A. Grepl and A.~T. Patera.
\newblock A posteriori error bounds for reduced-basis approximations of
  parametrized parabolic partial differential equations.
\newblock {\em ESAIM: Mathematical Modelling and Numerical Analysis},
  39(1):157--181, 2005.

\bibitem{Haasdonk2017Review}
B.~Haasdonk.
\newblock {\em Chapter 2: Reduced basis methods for parametrized {PDE}s--a
  tutorial introduction for stationary and instationary problems}, volume~15,
  pages 65--136.
\newblock SIAM Philadelphia, 2017.

\bibitem{HaasdonkOhlberger}
B.~Haasdonk and M.~Ohlberger.
\newblock {Reduced basis method for finite volume approximations of
  parametrized linear evolution equations}.
\newblock {\em M2AN Math. Model. Numer. Anal.}, 42(2):277--302, 2008.

\bibitem{HesthavenRozzaStammBook}
J.~S. Hesthaven, G.~Rozza, and B.~Stamm.
\newblock {\em Certified reduced basis methods for parametrized partial
  differential equations}.
\newblock SpringerBriefs in Mathematics. Springer, Cham; BCAM Basque Center for
  Applied Mathematics, Bilbao, 2016.
\newblock BCAM SpringerBriefs.

\bibitem{HKCHP}
D.~B.~P. Huynh, D.~J. Knezevic, Y.~Chen, J.~S. Hesthaven, and A.~T. Patera.
\newblock A natural-norm {Successive Constraint Method} for inf-sup lower
  bounds.
\newblock {\em CMAME}, 199:1963--1975, 2010.

\bibitem{HuynhSCM}
D.~B.~P. Huynh, G.~Rozza, S.~Sen, and A.~T. Patera.
\newblock {A successive constraint linear optimization method for lower bounds
  of parametric coercivity and inf-sup stability constants}.
\newblock {\em C. R. Acad. Sci. Paris, S$\acute{\rm e}$rie I.}, 345:473--478,
  2007.

\bibitem{JCX2018}
L.~Ji, Y.~Chen, and Z.~Xu.
\newblock A reduced basis method for the nonlinear {P}oisson-{B}oltzmann
  equation.
\newblock {\em Adv. Appl. Math. Mech}, 11(5):1200--1218, 2019.

\bibitem{Kunisch_Volkwein_POD}
K.~Kunisch and S.~Volkwein.
\newblock {Galerkin proper orthogonal decomposition methods for a general
  equation in fluid dynamics}.
\newblock {\em SIAM J. Numer. Anal.}, 40(2):492--515, 2002.

\bibitem{LiangPOD}
Y.~Liang, H.~Lee, S.~Lim, W.~Lin, K.~Lee, and C.~Wu.
\newblock Proper orthogonal decomposition and its applications-{Part I:
  Theory}.
\newblock {\em Journal of Sound and vibration}, 252(3):527--544, 2002.

\bibitem{negri2015efficient}
F.~Negri, A.~Manzoni, and D.~Amsallem.
\newblock Efficient model reduction of parametrized systems by matrix discrete
  empirical interpolation.
\newblock {\em J. Comput. Phys.}, 303:431--454, 2015.

\bibitem{nguyen2008efficient}
N.~C. Nguyen and J.~Peraire.
\newblock An efficient reduced-order modeling approach for non-linear
  parametrized partial differential equations.
\newblock {\em Int J Numer Methods Eng.}, 76(1):27--55, 2008.

\bibitem{nguyen2009reduced}
N.-C. Nguyen, G.~Rozza, and A.~T. Patera.
\newblock Reduced basis approximation and a posteriori error estimation for the
  time-dependent viscous {B}urgers' equation.
\newblock {\em Calcolo}, 46(3):157--185, 2009.

\bibitem{noor1979reduced}
A.~K. Noor and J.~M. Peters.
\newblock Reduced basis technique for nonlinear analysis of structures.
\newblock {\em AIAA J.}, 18(4):455--462, 1980.

\bibitem{peherstorfer2019sampling}
B.~Peherstorfer.
\newblock Sampling low-dimensional markovian dynamics for pre-asymptotically
  recovering reduced models from data with operator inference.
\newblock {\em arXiv preprint arXiv:1908.11233}, 2019.

\bibitem{PeherstorferButnaruWillcoxBungartz2014}
B.~Peherstorfer, D.~Butnaru, K.~Willcox, and H.~Bungartz.
\newblock Localized discrete empirical interpolation method.
\newblock {\em SIAM J. Sci. Comput.}, 36(1):A168--A192, 2014.

\bibitem{Quarteroni2015}
A.~Quarteroni, A.~Manzoni, and F.~Negri.
\newblock {\em Reduced basis methods for partial differential equations: {A}n
  introduction}, volume~92.
\newblock Springer, 2015.

\bibitem{Rozza2008}
G.~Rozza, D.~B.~P. Huynh, and A.~T. Patera.
\newblock Reduced basis approximation and a posteriori error estimation for
  affinely parametrized elliptic coercive partial differential equations.
\newblock {\em Arch. Comput. Methods Eng.}, 15(3):229--275, 2008.

\bibitem{ryckelynck2005priori}
D.~Ryckelynck.
\newblock A priori hyperreduction method: an adaptive approach.
\newblock {\em J. Comput. Phys.}, 202(1):346--366, 2005.

\bibitem{Ryckelynck2009}
D.~Ryckelynck.
\newblock Hyper-reduction of mechanical models involving internal variables.
\newblock {\em Int J Numer Methods Eng.}, 77(1):75--89, 2009.

\bibitem{veroy2003reduced}
K.~Veroy, C.~Prud'Homme, and A.~T. Patera.
\newblock Reduced-basis approximation of the viscous {B}urgers equation:
  rigorous a posteriori error bounds.
\newblock {\em C. R. Math.}, 337(9):619--624, 2003.

\bibitem{WillcoxPeraire2002}
K.~Willcox and J.~Peraire.
\newblock Balanced model reduction via the proper orthogonal decomposition.
\newblock {\em AIAA J.}, 40(11):2323--2330, 2002.

\end{thebibliography}

\end{document}